\documentclass[journal]{IEEEtran}

\usepackage[utf8]{inputenc}
\usepackage{color}
\usepackage{array}
\usepackage{verbatim}
\usepackage{float}
\usepackage{amsmath}
\usepackage{amsthm}
\usepackage{amssymb}
\usepackage{graphicx}
\usepackage{longtable}
\usepackage[unicode=true,
 bookmarks=false,
 breaklinks=false,pdfborder={0 0 1},colorlinks=false]
 {hyperref}
\hypersetup{
 colorlinks,bookmarksopen,bookmarksnumbered,citecolor=blue,urlcolor=blue}
\usepackage{cite}
\makeatletter
 \let\oldforeign@language\foreign@language
 \DeclareRobustCommand{\foreign@language}[1]{%
   \lowercase{\oldforeign@language{#1}}}

 \let\oldforeign@language\foreign@language
 \DeclareRobustCommand{\foreign@language}[1]{%
   \lowercase{\oldforeign@language{#1}}}

\ifCLASSINFOpdf
\else
\fi

\newcommand{\MYfooter}{\smash{
		\hfil\parbox[t][\height][t]{\textwidth}{\centering
			\thepage}\hfil\hbox{}}}

\makeatletter

\def\ps@IEEEtitlepagestyle{%
	\def\@oddhead{\parbox[t][\height][t]{\textwidth}{\centering \scriptsize
			Personal use of this material is permitted. Permission from the author(s) and/or copyright holder(s), must be obtained for all other uses. Please contact us and provide details if you believe this document breaches copyrights.\\
			\noindent\makebox[\linewidth]{}
		}\hfil\hbox{}}%
	\def\@evenhead{\scriptsize\thepage \hfil \leftmark\mbox{}}%
	\def\@oddfoot{\parbox[t][\height][l]{\textwidth}{
			\vspace{-20pt}{\rule{\textwidth}{0.4pt}}\\ \footnotesize\underline{To cite this article:}
			{\bf{\textcolor{red}{H. A. Hashim, L. J. Brown, and K. McIsaac, "Nonlinear Stochastic Attitude Filters on the Special Orthogonal Group 3: Ito and Stratonovich," IEEE Transactions on Systems, Man, and Cybernetics: Systems, vol. 49, no. 9, pp. 1853–1865, 2019.}}} doi: \href{https://doi.org/10.1109/TSMC.2018.2870290}{10.1109/TSMC.2018.2870290}\\
			\noindent\makebox[\linewidth]
		}\hfil\hbox{}}%
	\def\@evenfoot{\MYfooter}}

\makeatother
\newtheorem{defn}{Definition}

\newtheorem{lem}{Lemma}

\newtheorem{thm}{Theorem}

\newtheorem{assum}{Assumption}

\begin{document}
	\bstctlcite{IEEEexample:BSTcontrol}

\twocolumn
\title{Nonlinear Stochastic Attitude Filters on the \\
 Special Orthogonal Group 3: Ito and Stratonovich }

\author{Hashim~A.~Hashim, Lyndon J. Brown, and~Kenneth McIsaac
\thanks{H. A. Hashim, L. J. Brown and K. McIsaac are with the Department of Electrical and Computer Engineering,
University of Western Ontario, London, ON, Canada, N6A-5B9, e-mail: hmoham33@uwo.ca, lbrown@uwo.ca and kmcisaac@uwo.ca.}
}


\markboth{--,~Vol.~-, No.~-, \today}{Hashim \MakeLowercase{\textit{et al.}}: Nonlinear Stochastic Attitude Filter on the Special Orthogonal Group}
\markboth{}{Hashim \MakeLowercase{\textit{et al.}}: Nonlinear Stochastic Attitude Filter on the Special Orthogonal Group}

\maketitle

\begin{abstract}
This paper formulates the attitude filtering problem as a nonlinear
stochastic filter problem evolved directly on the Special Orthogonal
Group $\mathbb{SO}\left(3\right)$. One of the traditional potential
functions for nonlinear deterministic complimentary filters is studied
and examined against angular velocity measurements corrupted with
noise. This work demonstrates that the careful selection of the attitude
potential function allows to attenuate the noise associated with the
angular velocity measurements and results into superior convergence
properties of estimator and correction factor. The problem is formulated
as a stochastic problem through mapping $\mathbb{SO}\left(3\right)$
to Rodriguez vector parameterization. Two nonlinear stochastic complimentary
filters are developed on $\mathbb{SO}\left(3\right)$. The first stochastic
filter is driven in the sense of Ito and the second one considers
Stratonovich. The two proposed filters guarantee that errors in the
Rodriguez vector and estimates are semi-globally uniformly ultimately
bounded in mean square, and they converge to a small neighborhood
of the origin. Quaternion representation of the proposed observers is given. Simulation results are presented to illustrate the
effectiveness of the proposed filters considering high level of uncertainties
in angular velocity as well as body-frame vector measurements.
\end{abstract}


\IEEEpeerreviewmaketitle{}

\section{Introduction}

%
%
%
%
\IEEEPARstart{T}{his} paper concerns the problem of attitude estimation
of a rigid-body rotating in 3D space. In fact, attitude estimation
is one of the major sub-tasks in the field of robotics. The attitude
can be constructed from a set of vector measurements made on body-frame
and reference-frame as it acts as a linear transformation of one frame
to the other. Generally, the attitude estimation problem aims to minimize
the cost function such as Wahba’s Problem \cite{wahba1965least}.
The earliest work in \cite{wahba1965least} was purely algebraic.
Several alternative methods attempted to reconstruct the attitude
simply and statically by solving a set of simultaneous known inertial
and body-frame measurements, for instance, TRIAD or QUEST algorithms
\cite{black1964passive,shuster1981three} and singular value decomposition
(SVD) \cite{markley1988attitude}. However, vectorial measurements
are subject to significant noise and bias components. Therefore, the
category of static estimation in \cite{black1964passive,shuster1981three,markley1988attitude}
gives poor results in this case. Consequently, the attitude estimation
problem can be tackled either by Gaussian filter or nonlinear deterministic
filter.

In the last few decades, several Gaussian filters have been developed
mainly to obtain higher estimation performance with noise reduction.
Many attitude estimation algorithms are based on optimal stochastic
filtering for linear systems known as Kalman filter (KF) \cite{kalman1960new}.
The linearized version of KF can be modified in a certain way for
nonlinear systems to obtain the extended Kalman filter (EKF) \cite{anderson1979optimal}.
An early survey of attitude observers was presented in \cite{lefferts1982kalman}
and a more recent overview on attitude estimation was introduced in
\cite{crassidis2007survey}. Over the last three decades, several
nonlinear filters have been proposed to estimate the attitude of spacecrafts.
However, EKF and especially the multiplicative extended Kalman filter
(MEKF) is highly recommended and considered as a standard in most
spacecraft applications \cite{crassidis2007survey}. Generally, the
covariance of any noise components introduced in angular velocity
measurements is taken into account during filter design. The family
of KFs parameterize the global attitude problem using unit-quaternion.
The unit-quaternion provides a nonsingular attitude parameterization
of attitude matrix \cite{shuster1993survey}. Also, the unit-quaternion
kinematics and measurement models of the attitude can be defined by
a linear set of equations dependent on the quaternion state through
EKF. This advantage motivated researchers to employ the unit-quaternion
in attitude representation (for example \cite{lefferts1982kalman,markley2003attitude}).
Although EKF is subject to theoretical and practical problems, the
estimated state vector with the approximated covariance matrix gives
a reasonable estimate of uncertainties in the dynamics. In general,
a four-dimensional vector is used to describe a three-dimensional
one. Since, the covariance matrix associated with the quaternion vector
is $4\times4$, whereas the noise vector is $3\times1$, the covariance
is assumed to have rank 3. Generally, the state vector is $7\times1$
as it includes the four quaternion elements and the three bias components.
One of the earliest detailed derivations of EKF attitude design was
presented in \cite{lefferts1982kalman}. However, the unit-quaternion
kinematics and measurement models can be modified to suit KF with
a linear set of equations \cite{choukroun2006novel}. The KF in \cite{choukroun2006novel}
has the same state dimensions as EKF and to some degree, it can outperform
the EKF. MEKF \cite{markley2003attitude} is the modified version
of EKF and it is highly recommended for spacecraft applications \cite{crassidis2007survey}.
In MEKF, the true attitude state is the product of reference and estimated
error quaternion. The estimated error in quaternion is parameterized
from a three-dimensional vector in the body-frame, and the error is
estimated using EKF. Next, the MEKF is used to multiply the estimated
error and the reference quaternion. The estimated error should be
selected in such a way that it yields identity when multiplied by
the reference quaternion. The EKF can be modified into invariant extended
Kalman filter (IEKF), which has two groups of operations. The right
IEKF considers the errors modeled in the inertial-frame and the left
IEKF matches with the MEKF \cite{bonnable2009invariant}. IEKF has
autonomous error and its evolution error does not depend on the system
trajectory. A recently proposed attitude filtering solution known
as geometric approximate minimum-energy filter (GAMEF) approach \cite{zamani2013minimum}
is based on Mortensen’s deterministic minimum-energy \cite{mortensen1968maximum}.
Unlike KF, EKF, IEKF, and MEKF, the GAMEF kinematics are driven directly
on $\mathbb{SO}\left(3\right)$. In addition, KF, EKF, and IEKF are
based on first order optimal minimum-energy which makes them simpler
in computation and implementation. In contrast, MEKF and GAMEF are
second order optimal minimum-energy, and therefore they require more
calculation steps and more computational power. The Unscented Kalman
filter (UKF) uses the unit-quaternion kinematics, and its structure
is nearly similar to KF, however, UKF utilizes a set of sigma points
to enhance the probability distribution \cite{crassidis2003unscented}.
In spite of the fact that UKF requires less theoretical knowledge
and outperforms EKF in simulations, it requires more computational
power, while the sigma points could add complexity to the implementation
process \cite{haykin2001kalman}. Particle filters (PFs) belong to
the family of stochastic filters, but they do not follow the Gaussian
assumptions \cite{Arulampalam2002Particle}. The main idea of PFs
is the use of Monte-Carlo simulations for the weighted particle approximation
of the nonlinear distribution. In fact, PFs outperform EKF, however,
they have higher computational cost, and they do not fit small scale
systems \cite{crassidis2007survey}. Moreover, they do not have a
clear measure of how close the solution is to the optimal one \cite{zamani2013minimum}.
Quaternion based attitude PF showed a better performance than UKF
with higher processing calculations \cite{cheng2004particle}. All
the Gaussian filters described above as well as PFs are based on unit-quaternion,
where the main advantage is non-singularity in attitude parameterization,
while the main drawback is non-uniqueness in representation.

Aside from Gaussian filtering methods, nonlinear deterministic filters
provide an alternative solution of attitude estimation which aims
to establish convergence bounds with stable performance. Indeed, inertial
measurement units (IMUs) have a prominent role in enriching the research
of attitude estimation \cite{mahony2008nonlinear,rehbinder2004drift,metni2006attitude}.
IMUs fostered researchers to propose nonlinear deterministic complementary
filters on $\mathbb{SO}\left(3\right)$ using vectorial measurements
with the need of attitude reconstruction \cite{mahony2008nonlinear,lee2012exponential}
or directly from vectorial measurements without attitude reconstruction
\cite{mahony2008nonlinear,zlotnik2017nonlinear}. Also, the work done
in \cite{mahony2008nonlinear} provides the filter kinematics in quaternion
representation. 
In general, nonlinear deterministic filters achieve almost global
asymptotic stability as they disregard the noise impact in filter
derivation.

Nonlinear deterministic attitude filters have three distinctive advantages,
such as better tracking performance, less computational power, and
simplicity in derivation when compared to Gaussian filters or PFs
\cite{crassidis2007survey}. Furthermore, no sensor knowledge is required
in nonlinear deterministic filters, due to the fact that they omit
the noise component in filter derivation. Overall, nonlinear deterministic
attitude filters outperform Gaussian filters \cite{mahony2008nonlinear}.
Observers play a crucial role in different control applications, especially
for nonlinear stochastic systems (for example \cite{tong2011observer,zhou2017prescribed,dong2017observer}).
Aside from attitude observers, control applications are utilized for
nonlinear systems with uncertain components \cite{mohamed2014improved,hashim2017optimal}.
These applications could include robust stabilization \cite{wang2016approximate},
control of uncertain nonlinear multi-agent systems \cite{hashim2017neuro,hashim2017adaptive},
and stochastic nonlinear control for time-delay systems \cite{wang2017adaptive}.

Two major challenges have to be taken into account when designing
the attitude estimator, 1) the attitude problem of the rigid-body,
modeled on the Lie group of $\mathbb{SO}\left(3\right)$, is naturally
nonlinear; and 2) the true attitude kinematics rely on angular velocity.
Therefore, successful attitude estimation can be achieved by nonlinear
filter design relying on angular velocity measurements which are normally
contaminated with noise and bias components. Likewise, it is essential
that the estimator design considers any noise and/or bias components
introduced during the measurement process. Furthermore, any noise
component is characterized by randomness and irregular behavior. Having
this in mind, one of the traditional potential functions of nonlinear
deterministic complimentary filters evolved on $\mathbb{SO}\left(3\right)$
is studied (for example \cite{crassidis2007survey,mahony2008nonlinear})
taking into consideration angular velocity measurements corrupted
with bias and noise components. This study established that selecting
the potential function in an alternative way allowed to diminish the
noise. Hence, two nonlinear stochastic complementary filters on $\mathbb{SO}\left(3\right)$
are proposed here to improve the overall estimation quality. The first
stochastic filter is driven in the sense of Ito \cite{ito1984lectures}
and the second one is developed in the sense of Stratonovich \cite{stratonovich1967topics}.
In case when angular velocity measurement is contaminated with noise,
as far as the Rodriquez vector/($\mathbb{SO}\left(3\right)$) is concerned,
the proposed filters are able to 1) steer the error vector towards
an arbitrarily small neighborhood of the origin/(identity) in probability;
2) attenuate the noise impact to a very low level for known or unknown
bounded covariance; and 3) make the error semi-globally/(almost semi-globally)
uniformly ultimately bounded in mean square.

The rest of the paper is organized as follows: Section \ref{sec:SO3STCH_Math-Notations}
presents an overview of mathematical notation, $\mathbb{SO}\left(3\right)$
to Rodriguez vector parameterization, and some helpful properties
of the nonlinear stochastic attitude filter design. Attitude estimation
dynamic problem in Rodriguez vector with Gaussian noise vector which
satisfies the Brownian motion process is formulated in Section \ref{sec:SO3STCH_Problem-Formulation-in}.
The nonlinear stochastic filters on $\mathbb{SO}\left(3\right)$ and
the stability analysis are presented in Section \ref{sec:SO3STCH_Stochastic-Complementary-Filters}.
Section \ref{sec:SO3STCH_Simulations} shows the output performance
and discusses the simulation results of the proposed filters. Finally,
Section \ref{sec:SO3STCH_Conclusion} draws a conclusion of this work.

\section{Mathematical Notation \label{sec:SO3STCH_Math-Notations}}

Throughout this paper, $\mathbb{R}_{+}$ denotes the set of nonnegative
real numbers. $\mathbb{R}^{n}$ is the real $n$-dimensional space
while $\mathbb{R}^{n\times m}$ denotes the real $n\times m$ dimensional
space. For $x\in\mathbb{R}^{n}$, the Euclidean norm is defined as
$\left\Vert x\right\Vert =\sqrt{x^{\top}x}$, where $^{\top}$ is
the transpose of the associated component. $\mathcal{C}^{n}$ denotes
the set of functions with continuous $n$th partial derivatives. $\mathcal{K}$
denotes a set of continuous and strictly increasing functions such
that $\gamma:\mathbb{R}_{+}\rightarrow\mathbb{R}_{+}$ and vanishes
only at zero. $\mathcal{K}_{\infty}$ denotes a set of continuous
and strictly increasing functions which belongs to class $\mathcal{K}$
and is unbounded. $\mathbb{P}\left\{ \cdot\right\} $ denotes probability,
$\mathbb{E}\left[\cdot\right]$ denotes an expected value, and ${\rm exp}\left(\cdot\right)$
refers to an exponential of associated component. $\lambda\left(\cdot\right)$
is the set of singular values of the associated matrix with $\underline{\lambda}\left(\cdot\right)$
being the minimum value. $\mathbf{I}_{n}$ denotes identity matrix
with dimension $n$-by-$n$, and $\underline{\mathbf{0}}_{n}$ is
a zero vector with $n$-rows and one column. $V$ denotes a potential
function and for any $V\left(x\right)$ we have $V_{x}=\partial V/\partial x$
and $V_{xx}=\partial^{2}V/\partial x^{2}$.

Let $\mathbb{GL}\left(3\right)$ denote the 3 dimensional general
linear group which is a Lie group with smooth multiplication and inversion.
$\mathbb{SO}\left(3\right)$ denotes the Special Orthogonal Group
and is a subgroup of the general linear group. The attitude of a rigid-body
is defined as a rotational matrix $R$: 
\[
\mathbb{SO}\left(3\right)=\left\{ \left.R\in\mathbb{R}^{3\times3}\right|R^{\top}R=RR^{\top}=\mathbf{I}_{3}\text{, }{\rm det}\left(R\right)=1\right\} 
\]
where $\mathbb{I}_{n}$ is the identity matrix with $n$-dimensions
and ${\rm det\left(\cdot\right)}$ is the determinant of the associated
matrix. The associated Lie-algebra of $\mathbb{SO}\left(3\right)$
is termed $\mathfrak{so}\left(3\right)$ and is defined by 
\[
\mathfrak{so}\left(3\right)=\left\{ \left.\mathcal{A}\in\mathbb{R}^{3\times3}\right|\mathcal{A}^{\top}=-\mathcal{A}\right\} 
\]
with $\mathcal{A}$ being the space of skew-symmetric matrices and
define the map $\left[\cdot\right]_{\times}:\mathbb{R}^{3}\rightarrow\mathfrak{so}\left(3\right)$
such that 
\[
\mathcal{A}=\left[\alpha\right]_{\times}=\left[\begin{array}{ccc}
0 & -\alpha_{3} & \alpha_{2}\\
\alpha_{3} & 0 & -\alpha_{1}\\
-\alpha_{2} & \alpha_{1} & 0
\end{array}\right],\hspace{1em}\alpha=\left[\begin{array}{c}
\alpha_{1}\\
\alpha_{2}\\
\alpha_{3}
\end{array}\right]
\]
For all $\alpha,\beta\in\mathbb{R}^{3}$, we have $\left[\alpha\right]_{\times}\beta=\alpha\times\beta$
where $\times$ is the cross product between two vectors. Let the
vex operator be the inverse of $\left[\cdot\right]_{\times}$, denoted
by $\mathbf{vex}:\mathfrak{so}\left(3\right)\rightarrow\mathbb{R}^{3}$
such that $\mathbf{vex}\left(\mathcal{A}\right)=\alpha\in\mathbb{R}^{3}$.
Let $\boldsymbol{\mathcal{P}}_{a}$ denote the anti-symmetric projection
operator on the Lie-algebra $\mathfrak{so}\left(3\right)$, defined
by $\boldsymbol{\mathcal{P}}_{a}:\mathbb{R}^{3\times3}\rightarrow\mathfrak{so}\left(3\right)$
such that 
\[
\boldsymbol{\mathcal{P}}_{a}\left(\mathcal{B}\right)=\frac{1}{2}\left(\mathcal{B}-\mathcal{B}^{\top}\right)\in\mathfrak{so}\left(3\right)
\]
for all $\mathcal{B}\in\mathbb{R}^{3\times3}$. The following two
identities will be used in the subsequent derivations 
\begin{equation}
-\left[\beta\right]_{\times}\left[\alpha\right]_{\times}=\left(\beta^{\top}\alpha\right)\mathbf{I}_{3}-\alpha\beta^{\top},\quad\alpha,\beta\in{\rm \mathbb{R}}^{3}\label{eq:SO3STCH_Identity1}
\end{equation}
\begin{equation}
\left[R\alpha\right]_{\times}=R\left[\alpha\right]_{\times}R^{\top},\quad R\in\mathbb{SO}\left(3\right),\alpha\in\mathbb{R}^{3}\label{eq:SO3STCH_Identity2}
\end{equation}
The normalized Euclidean distance of a rotation matrix on $\mathbb{SO}\left(3\right)$
is given by the following equation 
\begin{equation}
\left\Vert R\right\Vert _{I}=\frac{1}{4}{\rm Tr}\left\{ \mathbf{I}_{3}-R\right\} \in\left[0,1\right]\label{eq:SO3STCH_Ecul_Dist}
\end{equation}
where ${\rm Tr}\left\{ \cdot\right\} $ denotes the trace of the associated
matrix and $\left\Vert R\right\Vert _{I}\in\left[0,1\right]$. The
attitude of a rigid-body can be constructed knowing angle of rotation
$\alpha\in\mathbb{R}$ and axis parameterization $u\in\mathbb{R}^{3}$.
This method of attitude reconstruction is termed angle-axis parameterization
\cite{shuster1993survey}. The mapping of angle-axis parameterization
to $\mathbb{SO}\left(3\right)$ is defined by $\mathcal{R}_{\alpha}:\mathbb{R}\times\mathbb{R}^{3}\rightarrow\mathbb{SO}\left(3\right)$
such that 
\begin{equation}
\mathcal{R}_{\alpha}\left(\alpha,u\right)=\mathbf{I}_{3}+\sin\left(\alpha\right)\left[u\right]_{\times}+\left(1-\cos\left(\alpha\right)\right)\left[u\right]_{\times}^{2}\label{eq:SO3STCH_att_ang}
\end{equation}
From the other side, the attitude can be defined knowing Rodriguez
parameters vector $\rho\in\mathbb{R}^{3}$. The associated map to
$\mathbb{SO}\left(3\right)$ is given by $\mathcal{R}_{\rho}:\mathbb{R}^{3}\rightarrow\mathbb{SO}\left(3\right)$
such that 
\begin{align}
\mathcal{R}_{\rho}\left(\rho\right)= & \frac{1}{1+\left\Vert \rho\right\Vert ^{2}}\left(\left(1-\left\Vert \rho\right\Vert ^{2}\right)\mathbf{I}_{3}+2\rho\rho^{\top}+2\left[\rho\right]_{\times}\right)\label{eq:SO3STCH_SO3_Rodr}
\end{align}
Substituting for the rotation matrix in \eqref{eq:SO3STCH_SO3_Rodr},
one can further show that the normalized Euclidean distance in \eqref{eq:SO3STCH_Ecul_Dist}
can be expressed in terms of Rodriguez parameters: 
\begin{equation}
\left\Vert R\right\Vert _{I}=\frac{1}{4}{\rm Tr}\left\{ \mathbf{I}_{3}-R\right\} =\frac{\left\Vert \rho\right\Vert ^{2}}{1+\left\Vert \rho\right\Vert ^{2}}\label{eq:SO3STCH_TR2}
\end{equation}
The anti-symmetric projection operator in square matrix space of the
rotation matrix $R$ in \eqref{eq:SO3STCH_SO3_Rodr} can be obtained
in the sense of Rodriguez parameters vector as 
\begin{align*}
\boldsymbol{\mathcal{P}}_{a}\left(R\right)= & 2\frac{1}{1+\left\Vert \rho\right\Vert ^{2}}\left[\rho\right]_{\times}\in\mathfrak{so}\left(3\right)
\end{align*}
It follows that the composition mapping $\boldsymbol{\Upsilon}_{a}\left(\cdot\right)$
is 
\begin{equation}
\boldsymbol{\Upsilon}_{a}\left(R\right)=\mathbf{vex}\left(\boldsymbol{\mathcal{P}}_{a}\left(R\right)\right)=2\frac{\rho}{1+\left\Vert \rho\right\Vert ^{2}}\in\mathbb{R}^{3}\label{eq:SO3STCH_VEX_Pa}
\end{equation}
where $\boldsymbol{\Upsilon}_{a}:=\mathbf{vex}\circ\boldsymbol{\mathcal{P}}_{a}$.

\section{Problem Formulation in Stochastic Sense \label{sec:SO3STCH_Problem-Formulation-in}}

\begin{figure}[h]
	\centering{}\includegraphics[scale=0.6]{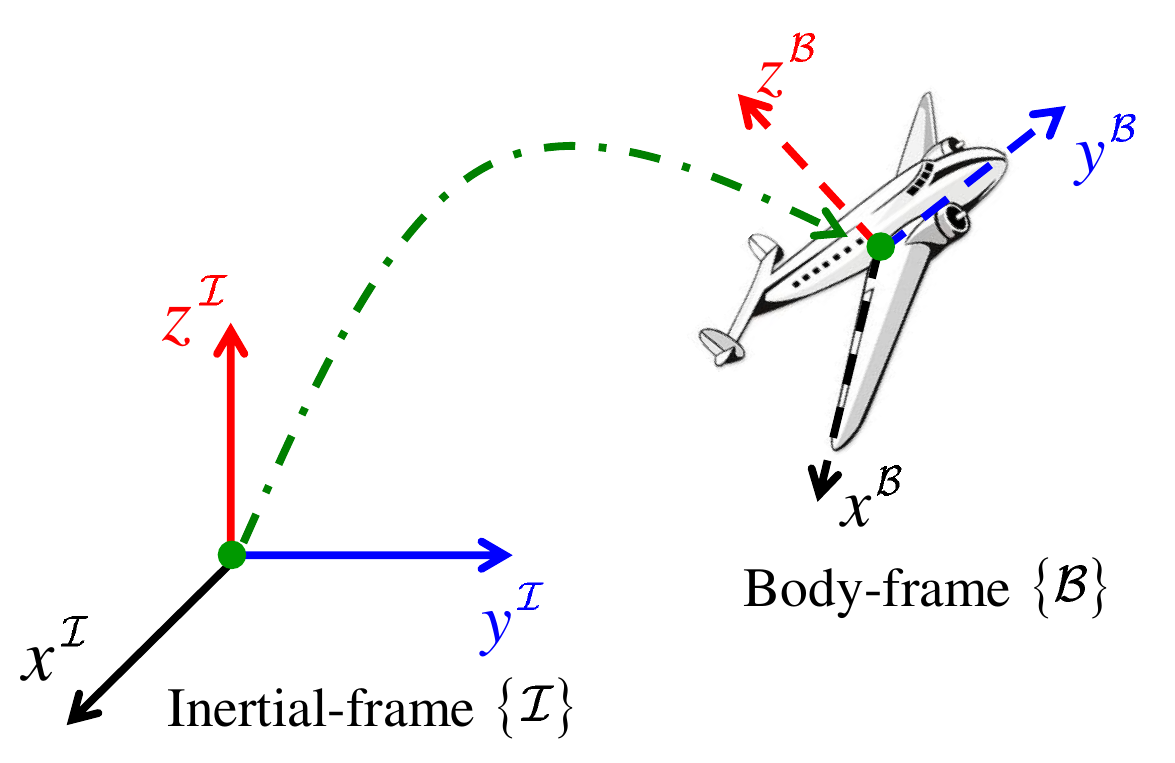}\caption{The orientation of a 3D rigid-body in body-frame relative to inertial-frame.}
	\label{fig:SO3STCH_0} 
\end{figure}

Let $R\in\mathbb{SO}\left(3\right)$ denote the attitude (rotational)
matrix, which describes the relative orientation of the body-frame
$\left\{ \mathcal{B}\right\} $ with respect to the inertial-frame
$\left\{ \mathcal{I}\right\} $ as given in Fig. \ref{fig:SO3STCH_0}.

The attitude can be extracted from $n$-known non-collinear inertial
vectors which are measured in a coordinate system fixed to the rigid
body. Let ${\rm v}_{i}^{\mathcal{B}}\in\mathbb{R}^{3}$ for $i=1,2,\ldots,n,$
be vectors measured in the body-fixed frame. Let $R\in\mathbb{SO}\left(3\right)$,
the body fixed-frame vector ${\rm v}_{i}^{\mathcal{B}}\in\mathbb{R}^{3}$
is defined by

\begin{equation}
{\rm v}_{i}^{\mathcal{B}}=R^{\top}{\rm v}_{i}^{\mathcal{I}}+{\rm b}_{i}^{\mathcal{B}}+\omega_{i}^{\mathcal{B}}\label{eq:SO3STCH_Vect_True}
\end{equation}
where ${\rm v}_{i}^{\mathcal{I}}\in\mathbb{R}^{3}$ denotes the inertial
fixed-frame vector for $i=1,2,\ldots,n$. ${\rm b}_{i}^{\mathcal{B}}$
and $\omega_{i}^{\mathcal{B}}$ denote the additive bias and noise
components of the associated body-frame vector, respectively, for
all ${\rm b}_{i}^{\mathcal{B}},\omega_{i}^{\mathcal{B}}\in\mathbb{R}^{3}$.
The assumption that $n\ge2$ is necessary for instantaneous three-dimensional
attitude determination. In case when $n=2$, the cross product of
the two measured vectors can be accounted as the third vector measurement
such that ${\rm v}_{3}^{\mathcal{I}}={\rm v}_{1}^{\mathcal{I}}\times{\rm v}_{2}^{\mathcal{I}}$
and ${\rm v}_{3}^{\mathcal{B}}={\rm v}_{1}^{\mathcal{B}}\times{\rm v}_{2}^{\mathcal{B}}$.
It is common to employ the normalized values of inertial and body-frame
vectors in the process of attitude estimation such as

\noindent 
\begin{equation}
\upsilon_{i}^{\mathcal{I}}=\frac{{\rm v}_{i}^{\mathcal{I}}}{\left\Vert {\rm v}_{i}^{\mathcal{I}}\right\Vert },\hspace{1em}\upsilon_{i}^{\mathcal{B}}=\frac{{\rm v}_{i}^{\mathcal{B}}}{\left\Vert {\rm v}_{i}^{\mathcal{B}}\right\Vert }\label{eq:SO3STCH_Vector_norm}
\end{equation}
In this manner, the attitude can be defined knowing $\upsilon_{i}^{\mathcal{I}}$
and $\upsilon_{i}^{\mathcal{B}}$. Gyroscope or the rate gyros measures
the angular velocity vector in the body-frame relative to the inertial-frame.
The measurement vector of angular velocity $\Omega_{m}\in\mathbb{R}^{3}$
is

\noindent 
\begin{equation}
\Omega_{m}=\Omega+b+\omega\label{eq:SO3STCH_Angular}
\end{equation}
where $\Omega\in\mathbb{R}^{3}$ denotes the true value of angular
velocity, $b$ denotes an unknown constant (bias) or slowly time-varying
vector, while $\omega$ denotes the noise component associated with
angular velocity measurements, for all $b,\omega\in\mathbb{R}^{3}$.
The noise vector $\omega$ is assumed to be Gaussian. The true attitude
dynamics and the associated Rodriguez vector dynamics are given in
\eqref{eq:SO3STCH_R_dynam} and \eqref{eq:SO3STCH_Rod_dynam}, respectively,
as 
\begin{equation}
\dot{R}=R\left[\Omega\right]_{\times}\label{eq:SO3STCH_R_dynam}
\end{equation}
\begin{equation}
\dot{\rho}=\frac{1}{2}\left(\mathbf{I}_{3}+\left[\rho\right]_{\times}+\rho\rho^{\top}\right)\Omega\label{eq:SO3STCH_Rod_dynam}
\end{equation}
In general, the measurement of angular velocity vector is subject
to additive noise and bias components. These components are characterized
by randomness and unknown behavior. In view of the fact that any unknown
components in angular velocity measurements may impair the estimation
process of the true attitude dynamics in \eqref{eq:SO3STCH_R_dynam}
or \eqref{eq:SO3STCH_Rod_dynam}, it is necessary to assume that the
attitude dynamics are excited by a wide-band of random Gaussian noise
process with zero mean. Combining angular velocity measurement in
\eqref{eq:SO3STCH_Angular} and the attitude dynamics in \eqref{eq:SO3STCH_Rod_dynam},
the attitude dynamics can be expressed as follows 
\begin{align}
\dot{\rho} & =\frac{1}{2}\left(\mathbf{I}_{3}+\left[\rho\right]_{\times}+\rho\rho^{\top}\right)\left(\Omega_{m}-b-\omega\right)\label{eq:SO3STCH_Stoch_non1}
\end{align}
where $\omega\in\mathbb{R}^{3}$ is a bounded continuous Gaussian
random noise vector with zero mean. The fact that derivative of any
Gaussian process yields Gaussian process allows us to write the stochastic
attitude dynamics as a function of Brownian motion process vector
$d\beta/dt\in\mathbb{R}^{3}$ \cite{khasminskii1980stochastic,jazwinski2007stochastic}.
Let $\left\{ \omega,t\geq t_{0}\right\} $ be a vector process of
independent Brownian motion process such that 
\begin{equation}
\omega=\mathcal{Q}\frac{d\beta}{dt}\label{eq:SO3STCH_noise}
\end{equation}
where $\mathcal{Q}\in\mathbb{R}^{3\times3}$ is an unknown time-variant
matrix with only nonzero and nonnegative bounded components in the
diagonal. The covariance component associated with the noise $\omega$
can be defined by $\mathcal{Q}^{2}=\mathcal{Q}\mathcal{Q}^{\top}$.
The properties of Brownian motion process are defined as \cite{ito1984lectures,jazwinski2007stochastic,deng2001stabilization}
\[
\mathbb{P}\left\{ \beta\left(0\right)=0\right\} =1,\hspace{1em}\mathbb{E}\left[d\beta/dt\right]=0,\hspace{1em}\mathbb{E}\left[\beta\right]=0
\]
Let the attitude dynamics of Rodriguez vector in \eqref{eq:SO3STCH_Rod_dynam}
be defined in the sense of Ito \cite{ito1984lectures}. Considering
the attitude dynamics in \eqref{eq:SO3STCH_Stoch_non1} and substituting
$\omega$ by $\mathcal{Q}d\beta/dt$ as in \eqref{eq:SO3STCH_noise},
the stochastic differential equation of \eqref{eq:SO3STCH_Rod_dynam}
in view of \eqref{eq:SO3STCH_Stoch_non1} can be expressed by 
\begin{align}
d\rho= & f\left(\rho,b\right)dt+g\left(\rho\right)\mathcal{Q}d\beta\label{eq:SO3STCH_Stoch_Ito}
\end{align}
Similarly, the stochastic dynamics of \eqref{eq:SO3STCH_R_dynam}
become 
\begin{equation}
dR=R\left[\Omega_{m}-b\right]_{\times}dt-R\left[\mathcal{Q}d\beta\right]_{\times}\label{eq:SO3STCH_R_Dyn1}
\end{equation}
where $b$ was defined in \eqref{eq:SO3STCH_Angular}, $g\left(\rho\right):=-\frac{1}{2}\left(\mathbf{I}_{3}+\left[\rho\right]_{\times}+\rho\rho^{\top}\right)$
and $f\left(\rho,b\right):=-g\left(\rho\right)\left(\Omega_{m}-b\right)$
with $g:\mathbb{R}^{3}\rightarrow\mathbb{R}^{3\times3}$ and $f:\mathbb{R}^{3}\times\mathbb{R}^{3}\rightarrow\mathbb{R}^{3}$.
$g\left(\rho\right)$ is locally Lipschitz in $\rho$, and $f\left(\rho,b\right)$
is locally Lipschitz in $\rho$ and $b$. Accordingly, the dynamic
system in \eqref{eq:SO3STCH_Stoch_Ito} has a solution for $t\in\left[t_{0},T\right]\forall t_{0}\leq T<\infty$
in the mean square sense and for any $\rho\left(t\right)\in\mathbb{R}^{3}$
such that $t\neq t_{0}$, $\rho-\rho_{0}$ is independent of $\left\{ \beta\left(\tau\right),\tau\geq t\right\} ,\forall t\in\left[t_{0},T\right]$
(Theorem 4.5 \cite{jazwinski2007stochastic}). Now the aim is to achieve
adaptive stabilization of an unknown bias and unknown time-variant
covariance matrix. Let $\sigma$ be the upper bound of $\mathcal{Q}^{2}$
such that 
\begin{equation}
\sigma=\left[{\rm max}\left\{ \mathcal{Q}_{1,1}^{2}\right\} ,{\rm max}\left\{ \mathcal{Q}_{2,2}^{2}\right\} ,{\rm max}\left\{ \mathcal{Q}_{3,3}^{2}\right\} \right]^{\top}\in\mathbb{R}^{3}\label{eq:SO3STCH_g_factor}
\end{equation}
where ${\rm max}\left\{ \cdot\right\} $ is the maximum value of an
element. 
\begin{assum}
	\label{Assum:SO3STCH_1} (Uniform boundedness of unknown parameters
	$b$ and $\sigma$) Let the vector $b$ and the nonnegative vector
	$\sigma$ belong to a given compact set $\Delta$ where $b,\sigma\in\Delta\subset\mathbb{R}^{3}$,
	and $b$ and $\sigma$ are upper bounded by a scalar $\Gamma$ such
	that $\left\Vert \Delta\right\Vert \leq\Gamma<\infty$. 
\end{assum}
\begin{defn}
	\label{def:SO3STCH_2} Consider the stochastic differential system
	in \eqref{eq:SO3STCH_Stoch_Ito}. For a given function $V\left(\rho\right)\in\mathcal{C}^{2}$,
	the differential operator $\mathcal{L}V$ is given by 
	\[
	\mathcal{L}V\left(\rho\right)=V_{\rho}^{\top}f\left(\rho,b\right)+\frac{1}{2}{\rm Tr}\left\{ g\left(\rho\right)\mathcal{Q}^{2}g^{\top}\left(\rho\right)V_{\rho\rho}\right\} 
	\]
	such that $V_{\rho}=\partial V/\partial\rho$, and $V_{\rho\rho}=\partial^{2}V/\partial\rho^{2}$. 
\end{defn}
\begin{defn}
	\label{def:SO3STCH_1}\cite{ji2006adaptive} The trajectory $\rho$
	of the stochastic differential system in \eqref{eq:SO3STCH_Stoch_Ito}
	is said to be semi-globally uniformly ultimately bounded (SGUUB) if
	for some compact set $\Lambda\in\mathbb{R}^{3}$ and any $\rho_{0}=\rho\left(t_{0}\right)$,
	there exists a constant $\kappa>0$, and a time constant $T=T\left(\kappa,\rho_{0}\right)$
	such that $\mathbb{E}\left[\left\Vert \rho\right\Vert \right]<\kappa,\forall t>t_{0}+T$. 
\end{defn}
\begin{lem}
	\label{lem:SO3STCH_1} \cite{deng2001stabilization,deng1997stochastic}
	Let the dynamic system in \eqref{eq:SO3STCH_Stoch_Ito} be assigned
	a potential function $V\in\mathcal{C}^{2}$ such that $V:\mathbb{R}^{3}\rightarrow\mathbb{R}_{+}$,
	class $\mathcal{K}_{\infty}$ function $\bar{\alpha}_{1}\left(\cdot\right)$
	and $\bar{\alpha}_{2}\left(\cdot\right)$, constants $c_{1}>0$ and
	$c_{2}\geq0$ and a nonnegative function $\mathbf{Z}\left(\left\Vert \rho\right\Vert \right)$
	such that 
	\begin{equation}
	\bar{\alpha}_{1}\left(\left\Vert \rho\right\Vert \right)\leq V\left(\rho\right)\leq\bar{\alpha}_{2}\left(\left\Vert \rho\right\Vert \right)\label{eq:SO3STCH_Vfunction_Lyap}
	\end{equation}
	\begin{align}
	\mathcal{L}V\left(\rho\right)= & V_{\rho}^{\top}f\left(\rho,b\right)+\frac{1}{2}{\rm Tr}\left\{ g\left(\rho\right)\mathcal{Q}^{2}g^{\top}\left(\rho\right)V_{\rho\rho}\right\} \nonumber \\
	\leq & -c_{1}\mathbf{Z}\left(\left\Vert \rho\right\Vert \right)+c_{2}\label{eq:SO3STCH_dVfunction_Lyap}
	\end{align}
	then for $\rho_{0}\in\mathbb{R}^{3}$, there exists almost a unique
	strong solution on $\left[0,\infty\right)$ for the dynamic system
	in \eqref{eq:SO3STCH_Stoch_Ito}, the solution $\rho$ is bounded
	in probability such that 
	\begin{equation}
	\mathbb{E}\left[V\left(\rho\right)\right]\leq V\left(\rho_{0}\right){\rm exp}\left(-c_{1}t\right)+\frac{c_{2}}{c_{1}}\label{eq:SO3STCH_EVfunction_Lyap}
	\end{equation}
	Furthermore, if the inequality in \eqref{eq:SO3STCH_EVfunction_Lyap}
	holds, then $\rho$ in \eqref{eq:SO3STCH_Stoch_Ito} is SGUUB in the
	mean square. In addition, when $c_{2}=0$, $f\left(0,b\right)=\underline{\mathbf{0}}_{3}$,
	$g\left(0\right)=\mathbf{0}_{3\times3}$, and $\mathbf{Z}\left(\left\Vert \rho\right\Vert \right)$
	is continuous, the equilibrium point $\rho=0$ is globally asymptotically
	stable in probability and the solution of $\rho$ satisfies 
	\begin{align}
	\mathbb{P}\left\{ \underset{t\rightarrow\infty}{{\rm lim}}\mathbf{Z}\left(\left\Vert \rho\right\Vert \right)=0\right\}  & =1,\hspace{1em}\forall\rho_{0}\in\mathbb{R}^{3}\label{eq:SO3STCH_PdVfunction_Lyap}
	\end{align}
\end{lem}
The proof of this lemma and existence of a unique solution can be
found in \cite{deng2001stabilization}. For a rotation matrix $R\in\mathbb{SO}\left(3\right)$,
let us define $\mathcal{U}\subseteq\mathbb{SO}\left(3\right)$ by
$\mathcal{U}:=\left\{ \left.R\right|{\rm Tr}\left\{ R\right\} =-1,\boldsymbol{\mathcal{P}}_{a}\left(R\right)=0\right\} $.
We have $-1\leq{\rm Tr}\left\{ R\right\} \leq3$ such that the set
$\mathcal{U}$ is forward invariant and unstable for the dynamic system
in \eqref{eq:SO3STCH_R_dynam} which implies that $\rho=\infty$.
For almost any initial condition such that $R_{0}\notin\mathcal{U}$
or $\rho_{0}\in\mathbb{R}^{3}$, we have $-1<{\rm Tr}\left\{ R_{0}\right\} \leq3$
and the trajectory of $\rho$ is semi-globally uniformly ultimately
bounded in mean square. 
\begin{lem}
	\label{lem:SO3STCH_2} (Young’s inequality) Let $x$ and $y$ be $x,y\in\mathbb{R}^{n}$.
	Then, for any $c>1$ and $d>1$ satisfying $\left(c-1\right)\left(d-1\right)=1$
	with a small positive constant $\varepsilon$, the following holds
	\begin{align}
	x^{\top}y & \leq\left(1/c\right)\varepsilon^{c}\left\Vert x\right\Vert ^{c}+\left(1/d\right)\varepsilon^{-d}\left\Vert y\right\Vert ^{d}\label{eq:SO3STCH_lem_ineq}
	\end{align}
\end{lem}
In the next section, the presence of noise will be examined in light
of a traditional form of potential function. The concept of an alternate
potential function with specific characteristics able to attenuate
the noise behavior will be carefully elucidated.

\section{Stochastic Complementary Filters On $\mathbb{SO}\left(3\right)$
	\label{sec:SO3STCH_Stochastic-Complementary-Filters}}

The main goal of attitude estimation is to derive the attitude estimate
$\hat{R}\rightarrow R$. Let's define the error in attitude estimate
from the body-frame to estimator-frame by 
\begin{equation}
\tilde{R}=R^{\top}\hat{R}\label{eq:SO3STCH_R_error}
\end{equation}
Let $\hat{b}$ and $\hat{\sigma}$ be estimates of unknown parameters
$b$ and $\sigma$, respectively. Define the error in vector $b$
and $\sigma$ by 
\begin{align}
\tilde{b} & =b-\hat{b}\label{eq:SO3STCH_b_tilde}\\
\tilde{\sigma} & =\sigma-\hat{\sigma}\label{eq:SO3STCH_sigma_tilde_strat}
\end{align}
Thus, driving $\hat{R}\rightarrow R$ ensures that $\tilde{R}\rightarrow\mathbf{I}_{3}$
and $\tilde{\rho}\rightarrow\underline{\mathbf{0}}_{3}$ where $\tilde{\rho}$
is Rodriguez error vector associated with $\tilde{R}$. In this section,
two nonlinear stochastic complementary filters are developed on the
Special Orthogonal Group $\mathbb{SO}\left(3\right)$. These filters
in the sense of Rodriguez vector guarantee that the error vector is
SGUUB in mean square for the case of noise contamination of the angular
velocity measurements.

\subsection{Nonlinear Deterministic Attitude Filter\label{subsec:SO3STCH_Det_Non_Fil}}

In this subsection, we aim to study the behavior of nonlinear deterministic
filter on $\mathbb{SO}\left(3\right)$ with noise introduced in angular
velocity measurements. The attitude $R$ can be reconstructed through
a set of measurements in \eqref{eq:SO3STCH_Vector_norm} to obtain
$R_{y}$, for instance \cite{black1964passive,shuster1981three,markley1988attitude}.
$R_{y}$ is corrupted with noise and bias greatly increase the difference
between $R_{y}$ and the true $R$. The filter design aims to use
the angular velocity measurements and the given $R_{y}$ to obtain
good estimate of $R$. Consider the following filter design 
\begin{align}
\dot{\hat{R}} & =\hat{R}\left[\Omega_{m}-\hat{b}-W\right]_{\times},\hspace{1em}\hat{R}\left(0\right)=\hat{R}_{0}\label{eq:SO3STCH_dRest_Det}\\
\dot{\hat{b}} & =\gamma_{1}\boldsymbol{\Upsilon}_{a}(\tilde{R}),\hspace{1em}\hat{b}\left(0\right)=\hat{b}_{0},\tilde{R}=R_{y}^{\top}\hat{R}\label{eq:SO3STCH_best_Det}\\
W & =k_{1}\boldsymbol{\Upsilon}_{a}(\tilde{R}),\hspace{1em}\tilde{R}=R_{y}^{\top}\hat{R}\label{eq:SO3STCH_Wcorr_Det}
\end{align}
where $\Omega_{m}$ is angular velocity measurement, $\hat{b}\in\mathbb{R}^{3}$
is the estimate of the unknown bias $b$, and $\boldsymbol{\Upsilon}_{a}(\tilde{R})=\mathbf{vex}\left(\boldsymbol{\mathcal{P}}_{a}(\tilde{R})\right)$
was given in \eqref{eq:SO3STCH_VEX_Pa}. Also, $\gamma_{1}>0$ is
an adaptation gain and $k_{1}$ is a positive constant.

Let the error in vector $b$ be defined as in \eqref{eq:SO3STCH_b_tilde}
and assume that no noise was introduced to the dynamics $\left(\omega=\underline{\mathbf{0}}_{3}\right)$.
The derivative of attitude error in \eqref{eq:SO3STCH_R_error} can
be obtained from \eqref{eq:SO3STCH_R_dynam} and \eqref{eq:SO3STCH_dRest_Det}
as 
\begin{align}
\dot{\tilde{R}}= & \tilde{R}\left[\Omega-\tilde{R}^{\top}\Omega+\tilde{b}-W\right]_{\times}\label{eq:SO3STCH_dR_tilde_Det}
\end{align}
where $\left[\tilde{R}^{\top}\Omega\right]_{\times}=\tilde{R}^{\top}\left[\Omega\right]_{\times}\tilde{R}$.
Hence, in view of \eqref{eq:SO3STCH_R_Dyn1} and \eqref{eq:SO3STCH_Stoch_Ito},
the error dynamic in \eqref{eq:SO3STCH_dR_tilde_Det} can be expressed
in Rodriguez error vector dynamic by 
\begin{align}
\dot{\tilde{\rho}} & =\frac{1}{2}\left(\mathbf{I}_{3}+\left[\tilde{\rho}\right]_{\times}+\tilde{\rho}\tilde{\rho}^{\top}\right)\left(\Omega-\tilde{R}^{\top}\Omega+\tilde{b}-W\right)\label{eq:SO3STCH_dr_tilde_Det}
\end{align}
From literature, one of traditional potential functions for adaptive
filter estimation is $V\left(\tilde{R},\tilde{b}\right)=\frac{1}{4}{\rm Tr}\left\{ \mathbf{I}_{3}-\tilde{R}\right\} +\frac{1}{2\gamma_{1}}\tilde{b}^{\top}\tilde{b}$
(for example \cite{crassidis2007survey,mahony2008nonlinear}). The
equivalent of the aforementioned function in form of Rodriguez error
is 
\begin{align}
V\left(\tilde{\rho},\tilde{b}\right)= & \frac{\left\Vert \tilde{\rho}\right\Vert ^{2}}{1+\left\Vert \tilde{\rho}\right\Vert ^{2}}+\frac{1}{2\gamma_{1}}\tilde{b}^{\top}\tilde{b}\label{eq:SO3STCH_LyapV_Det}
\end{align}
let $\tilde{f}:=\frac{1}{2}\left(\mathbf{I}_{3}+\left[\tilde{\rho}\right]_{\times}+\tilde{\rho}\tilde{\rho}^{\top}\right)\left(\Omega-\tilde{R}^{\top}\Omega+\tilde{b}-W\right)$.
For $V:=V\left(\tilde{\rho},\tilde{b}\right)$, the derivative of
\eqref{eq:SO3STCH_LyapV_Det} is 
\begin{align}
\dot{V} & =V_{\tilde{\rho}}^{\top}\tilde{f}-\frac{1}{\gamma_{1}}\tilde{b}^{\top}\dot{\hat{b}}\label{eq:SO3STCH_Vdot_Det}\\
& =\boldsymbol{\Upsilon}_{a}(\tilde{R})^{\top}\left(\tilde{b}-W\right)-\frac{1}{\gamma_{1}}\tilde{b}^{\top}\dot{\hat{b}}\nonumber 
\end{align}
where $\frac{1}{2}V_{\mathcal{\tilde{\rho}}}^{\top}\left(\mathbf{I}_{3}+\left[\tilde{\rho}\right]_{\times}+\tilde{\rho}\tilde{\rho}^{\top}\right)\left(\Omega-\tilde{R}^{\top}\Omega\right)=0$
which was obtained by substitution of $\tilde{R}=\mathcal{R}_{\tilde{\rho}}\left(\tilde{\rho}\right)$
in \eqref{eq:SO3STCH_SO3_Rodr}. Substituting for $\dot{\hat{b}}$
and $W$ in \eqref{eq:SO3STCH_best_Det} and \eqref{eq:SO3STCH_Wcorr_Det},
respectively, yields 
\begin{align}
\dot{V} & =-k_{1}\left\Vert \boldsymbol{\Upsilon}_{a}(\tilde{R})\right\Vert ^{2}=-4k_{1}\frac{\left\Vert \tilde{\rho}\right\Vert ^{2}}{\left(1+\left\Vert \tilde{\rho}\right\Vert ^{2}\right)^{2}}\label{eq:SO3STCH_Vdot2_Det}
\end{align}
Lyapunov’s direct method ensures that for ${\rm Tr}\left\{ \tilde{R}_{0}\right\} \neq-1$,
$\boldsymbol{\Upsilon}_{a}(\tilde{R})$ converges asymptotically
to zero. As such, $\left(\mathbf{I}_{3},\underline{\mathbf{0}}_{3}\right)$
is an isolated equilibrium point and $\left(\tilde{R},\tilde{b}\right)\rightarrow\left(\mathbf{I}_{3},\underline{\mathbf{0}}_{3}\right)$
for $\omega=\underline{\mathbf{0}}_{3}$ \cite{mahony2008nonlinear}.
If angular velocity measurements ($\Omega_{m}$) are contaminated
with noise $\left(\omega\neq\underline{\mathbf{0}}_{3}\right)$, it
is more convenient to represent the differential operator in \eqref{eq:SO3STCH_Vdot_Det}
in the form of Definition \ref{def:SO3STCH_2}. Hence, the following
extra term will appear 
\[
\frac{1}{2}{\rm Tr}\left\{ \tilde{g}^{\top}V_{\tilde{\rho}\tilde{\rho}}\tilde{g}\mathcal{Q}^{2}\right\} =\frac{1}{4\left(1+\left\Vert \tilde{\rho}\right\Vert ^{2}\right)}{\rm Tr}\left\{ \left(\mathbf{I}_{3}-3\tilde{\rho}\tilde{\rho}^{\top}\right)\mathcal{Q}^{2}\right\} 
\]
In this case, the operator $\mathcal{L}V\left(0,0\right)=\frac{1}{4}{\rm Tr}\left\{ \mathcal{Q}^{2}\right\} $
which implies that the significant impact of covariance matrix $\mathcal{Q}^{2}$
cannot be lessened. One way to attenuate the noise associated with
the angular velocity measurements is to chose a potential function
in the sense of Rodriguez error vector $\tilde{\rho}$ of order higher
than two. It is worth mentioning that the deterministic filter in
\eqref{eq:SO3STCH_dRest_Det}, \eqref{eq:SO3STCH_best_Det} and \eqref{eq:SO3STCH_Wcorr_Det}
is known as a passive complementary filter proposed in \cite{mahony2008nonlinear}.

\subsection{Nonlinear Stochastic Attitude Filter in Ito Sense\label{subsec:SO3STCH_Nonlinear-Ito}}

Generally, the assumption behind nonlinear deterministic filters is
that angular velocity vector measurements are joined with constant
or slowly time-variant bias \cite{crassidis2007survey,mahony2008nonlinear}.
However, angular velocity vector measurements are typically subject
to additive noise components which may weaken the estimation process
of the true attitude dynamics in \eqref{eq:SO3STCH_R_dynam}. Therefore,
we aim to design a nonlinear stochastic filter in Ito sense taking
into consideration that angular velocity vector measurements are subject
to a constant bias and a wide-band of Gaussian random with zero mean
such that $\mathbb{E}\left[\omega\right]=0$. Let the true inertial
vector ${\rm v}_{i}^{\mathcal{I}}$ and body-frame vector ${\rm v}_{i}^{\mathcal{B}}$
be defined as in \eqref{eq:SO3STCH_Vect_True}. Let the error in attitude
estimate be similar to \eqref{eq:SO3STCH_R_error}.

Consider the nonlinear stochastic filter design 
\begin{align}
\dot{\hat{R}}= & \hat{R}\left[\Omega_{m}-\hat{b}-W\right]_{\times},\hspace{1em}\hat{R}\left(0\right)=\hat{R}_{0}\label{eq:SO3STCH_dRest_ito}\\
\dot{\hat{b}}= & \gamma_{1}||\tilde{R}||_{I}\boldsymbol{\Upsilon}_{a}(\tilde{R})-\gamma_{1}k_{b}\hat{b},\hspace{1em}\hat{b}\left(0\right)=\hat{b}_{0}\label{eq:SO3STCH_best_ito}\\
\mathcal{\dot{\hat{\sigma}}}= & k_{1}\gamma_{2}||\tilde{R}||_{I}\mathcal{D}_{\Upsilon}^{\top}\boldsymbol{\Upsilon}_{a}(\tilde{R})-\gamma_{2}k_{\sigma}\hat{\sigma},\hspace{1em}\hat{\sigma}\left(0\right)=\hat{\sigma}_{0}\label{eq:SO3STCH_sest_ito}\\
W= & \frac{k_{1}}{\varepsilon}\frac{2-||\tilde{R}||_{I}}{1-||\tilde{R}||_{I}}\boldsymbol{\Upsilon}_{a}(\tilde{R})+k_{2}\mathcal{D}_{\Upsilon}\hat{\sigma}\label{eq:SO3STCH_Wcorr_ito}
\end{align}
where $\Omega_{m}$ is angular velocity measurement defined in \eqref{eq:SO3STCH_Angular},
$\hat{b}$ is the estimate of the unknown bias $b$, $\hat{\sigma}$
is the estimate of $\sigma$ which includes the upper bound of $\mathcal{Q}^{2}$
as given in \eqref{eq:SO3STCH_g_factor}, $\tilde{R}=R_{y}^{\top}\hat{R}$
with $R_{y}$ being the reconstructed attitude, $\boldsymbol{\Upsilon}_{a}(\tilde{R})=\mathbf{vex}\left(\boldsymbol{\mathcal{P}}_{a}(\tilde{R})\right)$
as given in \eqref{eq:SO3STCH_VEX_Pa}, $\mathcal{D}_{\Upsilon}=\left[\boldsymbol{\Upsilon}_{a}(\tilde{R}),\boldsymbol{\Upsilon}_{a}(\tilde{R}),\boldsymbol{\Upsilon}_{a}(\tilde{R})\right]$,
and $||\tilde{R}||_{I}$ is the Euclidean distance of $\tilde{R}$
as defined in \eqref{eq:SO3STCH_Ecul_Dist}. Also, $\gamma_{1}>0$
and $\gamma_{2}>0$ are adaptation gains, $\varepsilon>0$ is a small
constant, while $k_{b}$, $k_{\sigma}$, $k_{1}$ and $k_{2}$ are
positive constants. Quaternion representation of Ito's filter is presented in \nameref{sec:SO3STCH_AppendixB}. 
\begin{thm}
	\textbf{\label{thm:SO3STCH_1} }Consider the rotation dynamics in
	\eqref{eq:SO3STCH_R_Dyn1}, angular velocity measurements in \eqref{eq:SO3STCH_Angular}
	in addition to other given vectorial measurements in \eqref{eq:SO3STCH_Vector_norm}
	coupled with the observer \eqref{eq:SO3STCH_dRest_ito}, \eqref{eq:SO3STCH_best_ito},
	\eqref{eq:SO3STCH_sest_ito}, and \eqref{eq:SO3STCH_Wcorr_ito}. Assume
	that two or more body-frame non-collinear vectors are available for
	measurements and the design parameters $\gamma_{1}$ , $\gamma_{2}$
	, $\varepsilon$, $k_{b}$, $k_{\sigma}$, $k_{1}$, and $k_{2}$
	are chosen appropriately with $\varepsilon$ being selected sufficiently
	small. Then, for angular velocity measurements contaminated with noise $\left(\omega\neq\underline{\mathbf{0}}_{3}\right)$,
	all the signals in the closed-loop system is semi-globally uniformly
	ultimately bounded in mean square. In addition, the observer errors
	can be minimized by the appropriate selection of the design parameters.
\end{thm}
\textbf{Proof: }Let the error in vector $b$ be defined as in \eqref{eq:SO3STCH_b_tilde}.
Therefore, the derivative of attitude error in incremental form of
\eqref{eq:SO3STCH_R_error} can be obtained from \eqref{eq:SO3STCH_Stoch_Ito}
and \eqref{eq:SO3STCH_dRest_ito} by 
\begin{align}
d\tilde{R}= & R^{\top}\hat{R}\left[\Omega_{m}-\hat{b}-W\right]_{\times}dt+\left[\Omega\right]_{\times}^{\top}R^{\top}\hat{R}dt\nonumber \\
= & \left(\tilde{R}\left[\Omega\right]_{\times}+\left[\Omega\right]_{\times}^{\top}\tilde{R}+\tilde{R}\left[\tilde{b}-W\right]_{\times}\right)dt+\tilde{R}\left[\mathcal{Q}d\beta\right]_{\times}\nonumber \\
= & \tilde{R}\left[\Omega-\tilde{R}^{\top}\Omega+\tilde{b}-W\right]_{\times}dt+\tilde{R}\left[\mathcal{Q}d\beta\right]_{\times}\label{eq:SO3STCH_dR_tilde_ito}
\end{align}
Similar extraction of Rodriguez error vector dynamic in view of \eqref{eq:SO3STCH_R_Dyn1}
to \eqref{eq:SO3STCH_Stoch_Ito} can be expressed from \eqref{eq:SO3STCH_dR_tilde_ito}
to \eqref{eq:SO3STCH_StochErr_ito} in Ito's representation \cite{ito1984lectures}
as

\begin{align}
d\tilde{\rho}= & \tilde{f}dt+\tilde{g}\mathcal{Q}d\beta\label{eq:SO3STCH_StochErr_ito}
\end{align}
where $\tilde{\rho}$ is the Rodriguez error vector associated with
$\tilde{R}$. Let $\tilde{g}=\frac{1}{2}\left(\mathbf{I}_{3}+\left[\tilde{\rho}\right]_{\times}+\tilde{\rho}\tilde{\rho}^{\top}\right)$
and $\tilde{f}=\tilde{g}\left(\Omega-\tilde{R}^{\top}\Omega+\tilde{b}-W\right)$.
Consider the following potential function 
\begin{align}
V\left(\tilde{\rho},\tilde{b},\tilde{\sigma}\right)= & \left(\frac{\left\Vert \tilde{\rho}\right\Vert ^{2}}{1+\left\Vert \tilde{\rho}\right\Vert ^{2}}\right)^{2}+\frac{1}{2\gamma_{1}}\tilde{b}^{\top}\tilde{b}+\frac{1}{2\gamma_{2}}\tilde{\sigma}^{\top}\tilde{\sigma}\label{eq:SO3STCH_LyapV_ito}
\end{align}
For $V:=V\left(\tilde{\rho},\tilde{b},\tilde{\sigma}\right)$, the
differential operator $\mathcal{L}V$ in Definition \ref{def:SO3STCH_2}
for the dynamic system in \eqref{eq:SO3STCH_StochErr_ito} can be
expressed as 
\begin{equation}
\mathcal{L}V=V_{\tilde{\rho}}^{\top}\tilde{f}+\frac{1}{2}{\rm Tr}\left\{ \tilde{g}^{\top}V_{\tilde{\rho}\tilde{\rho}}\tilde{g}\mathcal{Q}^{2}\right\} -\frac{1}{\gamma_{1}}\tilde{b}^{\top}\dot{\hat{b}}-\frac{1}{\gamma_{2}}\tilde{\sigma}^{\top}\dot{\hat{\sigma}}\label{eq:SO3STCH_LV_Function_ito}
\end{equation}
where $V_{\tilde{\rho}}=\partial V/\partial\tilde{\rho}$ and $V_{\tilde{\rho}\tilde{\rho}}=\partial V^{2}/\partial^{2}\tilde{\rho}$.
The first and the second partial derivatives of \eqref{eq:SO3STCH_LyapV_ito}
with respect to $\tilde{\rho}$ can be obtained as follows 
\begin{align}
V_{\mathcal{\tilde{\rho}}}= & 4\frac{\left\Vert \tilde{\rho}\right\Vert ^{2}}{\left(1+\left\Vert \tilde{\rho}\right\Vert ^{2}\right)^{3}}\tilde{\rho}\label{eq:SO3STCH_LyapVv_ito}\\
V_{\mathcal{\tilde{\rho}}\mathcal{\tilde{\rho}}}= & 4\frac{\left(1+\left\Vert \tilde{\rho}\right\Vert ^{2}\right)\left\Vert \tilde{\rho}\right\Vert ^{2}\mathbf{I}_{3}+\left(2-4\left\Vert \tilde{\rho}\right\Vert ^{2}\right)\tilde{\rho}\tilde{\rho}^{\top}}{\left(1+\left\Vert \tilde{\rho}\right\Vert ^{2}\right)^{4}}\label{eq:SO3STCH_LyapVvv_ito}
\end{align}
substituting $\tilde{R}=\mathcal{R}_{\tilde{\rho}}\left(\tilde{\rho}\right)$
in \eqref{eq:SO3STCH_SO3_Rodr}, one can verify that 
\[
\frac{1}{2}V_{\mathcal{\tilde{\rho}}}^{\top}\left(\mathbf{I}_{3}+\left[\tilde{\rho}\right]_{\times}+\tilde{\rho}\tilde{\rho}^{\top}\right)\left(\Omega-\tilde{R}^{\top}\Omega\right)=0
\]
Hence, the first part of the differential operator $\mathcal{L}V$
in \eqref{eq:SO3STCH_LV_Function_ito} can be evaluated by 
\begin{align}
V_{\mathcal{\tilde{\rho}}}^{\top}\tilde{f} & =2\frac{\left\Vert \tilde{\rho}\right\Vert ^{2}}{\left(1+\left\Vert \tilde{\rho}\right\Vert ^{2}\right)^{2}}\tilde{\rho}^{\top}\left(\tilde{b}-W\right)\label{eq:SO3STCH_LyapVdot_ito}
\end{align}
Keeping in mind the identity in \eqref{eq:SO3STCH_Identity1} and
$\tilde{g}$ in \eqref{eq:SO3STCH_StochErr_ito} and combining them
with \eqref{eq:SO3STCH_LyapVvv_ito}, the component ${\rm Tr}\left\{ \tilde{g}^{\top}V_{\tilde{\rho}\tilde{\rho}}\tilde{g}\mathcal{Q}^{2}\right\} $
can be simplified and expressed as 
\begin{align}
\frac{1}{2}{\rm Tr}\left\{ \tilde{g}^{\top}V_{\tilde{\rho}\tilde{\rho}}\tilde{g}\mathcal{Q}^{2}\right\} = & \frac{1}{2\left(1+\left\Vert \tilde{\rho}\right\Vert ^{2}\right)^{3}}{\rm Tr}\left\{ \left(1+\left\Vert \tilde{\rho}\right\Vert ^{2}\right)\left\Vert \tilde{\rho}\right\Vert ^{2}\mathcal{Q}^{2}\right.\nonumber \\
& \left.+\left(2-\left\Vert \tilde{\rho}\right\Vert ^{2}-3\left\Vert \tilde{\rho}\right\Vert ^{4}\right)\tilde{\rho}\tilde{\rho}^{\top}\mathcal{Q}^{2}\right\} \label{eq:SO3STCH_LyapTrace_ito}
\end{align}
Let $\bar{q}=\left[\mathcal{Q}_{1,1},\mathcal{Q}_{2,2},\mathcal{Q}_{3,3}\right]^{\top}$
and $\sigma$ be similar to \eqref{eq:SO3STCH_g_factor}. From \eqref{eq:SO3STCH_LyapVdot_ito}
and \eqref{eq:SO3STCH_LyapTrace_ito}, one can write the operator
$\mathcal{L}V$ in \eqref{eq:SO3STCH_LV_Function_ito} as{\small{}
	\begin{align}
	\mathcal{L}V= & 2\frac{\left\Vert \tilde{\rho}\right\Vert ^{2}\tilde{\rho}^{\top}\left(\tilde{b}-W\right)}{\left(1+\left\Vert \tilde{\rho}\right\Vert ^{2}\right)^{2}}+\frac{{\rm Tr}\left\{ \left(2-\left\Vert \tilde{\rho}\right\Vert ^{2}-3\left\Vert \tilde{\rho}\right\Vert ^{4}\right)\tilde{\rho}\tilde{\rho}^{\top}\mathcal{Q}^{2}\right\} }{2\left(1+\left\Vert \tilde{\rho}\right\Vert ^{2}\right)^{3}}\nonumber \\
	& +\frac{{\rm Tr}\left\{ \left\Vert \tilde{\rho}\right\Vert ^{2}\mathcal{Q}^{2}\right\} }{2\left(1+\left\Vert \tilde{\rho}\right\Vert ^{2}\right)^{2}}-\frac{1}{\gamma_{1}}\tilde{b}^{\top}\dot{\hat{b}}-\frac{1}{\gamma_{2}}\tilde{\sigma}^{\top}\dot{\hat{\sigma}}\label{eq:SO3STCH_LV_ito_Pre}
	\end{align}
}Since $\left\Vert \bar{q}\right\Vert ^{2}={\rm Tr}\left\{ \mathcal{Q}^{2}\right\} $
and ${\rm Tr}\left\{ \tilde{\rho}\tilde{\rho}^{\top}\mathcal{Q}^{2}\right\} \leq\left\Vert \tilde{\rho}\right\Vert ^{2}\left\Vert \bar{q}\right\Vert ^{2}$,
we have 
\begin{align}
& \mathcal{L}V\leq2\frac{\left\Vert \tilde{\rho}\right\Vert ^{2}\tilde{\rho}^{\top}\left(\tilde{b}-W\right)}{\left(1+\left\Vert \tilde{\rho}\right\Vert ^{2}\right)^{2}}+\frac{\left\Vert \tilde{\rho}\right\Vert ^{4}{\rm Tr}\left\{ \mathcal{Q}^{2}\right\} +3\left\Vert \tilde{\rho}\right\Vert ^{2}\left\Vert \bar{q}\right\Vert ^{2}}{2\left(1+\left\Vert \tilde{\rho}\right\Vert ^{2}\right)^{3}}\nonumber \\
& \hspace{1.5em}-\frac{\left(1+3\left\Vert \tilde{\rho}\right\Vert ^{2}\right)\left\Vert \tilde{\rho}\right\Vert ^{2}\tilde{\rho}^{\top}\mathcal{Q}^{2}\tilde{\rho}}{2\left(1+\left\Vert \tilde{\rho}\right\Vert ^{2}\right)^{3}}-\frac{1}{\gamma_{1}}\tilde{b}^{\top}\dot{\hat{b}}-\frac{1}{\gamma_{2}}\tilde{\sigma}^{\top}\dot{\hat{\sigma}}\label{eq:SO3STCH_LV_ito}
\end{align}
According to Lemma \ref{lem:SO3STCH_2}, the following equation holds
\begin{align}
\frac{3\left\Vert \tilde{\rho}\right\Vert ^{2}\left\Vert \bar{q}\right\Vert ^{2}}{2\left(1+\left\Vert \tilde{\rho}\right\Vert ^{2}\right)^{3}} & \leq\frac{9\left\Vert \tilde{\rho}\right\Vert ^{4}}{8\left(1+\left\Vert \tilde{\rho}\right\Vert ^{2}\right)^{6}\varepsilon}+\frac{\varepsilon}{2}\left\Vert \bar{q}\right\Vert ^{4}\nonumber \\
& \leq\frac{9\left\Vert \tilde{\rho}\right\Vert ^{4}}{8\left(1+\left\Vert \tilde{\rho}\right\Vert ^{2}\right)^{3}\varepsilon}+\frac{\varepsilon}{2}\left(\sum_{i=1}^{3}\sigma_{i}\right)^{2}\label{eq:SO3STCH_Inq_Lyap}
\end{align}
where $\varepsilon$ is a sufficiently small positive constant. Combining
\eqref{eq:SO3STCH_Inq_Lyap} with \eqref{eq:SO3STCH_LV_ito} yields
\begin{align}
\mathcal{L}V\leq & 2\frac{\left\Vert \tilde{\rho}\right\Vert ^{2}\tilde{\rho}^{\top}\left(\tilde{b}-W\right)}{\left(1+\left\Vert \tilde{\rho}\right\Vert ^{2}\right)^{2}}+2\frac{\frac{1}{4}\sum_{i=1}^{3}\sigma_{i}+\frac{9}{16\varepsilon}}{\left(1+\left\Vert \tilde{\rho}\right\Vert ^{2}\right)^{3}}\left\Vert \tilde{\rho}\right\Vert ^{4}\nonumber \\
& -\frac{1}{\gamma_{1}}\tilde{b}^{\top}\dot{\hat{b}}-\frac{1}{\gamma_{2}}\tilde{\sigma}^{\top}\dot{\hat{\sigma}}\nonumber \\
& -\frac{\left(1+3\left\Vert \tilde{\rho}\right\Vert ^{2}\right)\left\Vert \tilde{\rho}\right\Vert ^{2}\tilde{\rho}^{\top}\mathcal{Q}^{2}\tilde{\rho}}{2\left(1+\left\Vert \tilde{\rho}\right\Vert ^{2}\right)^{3}}+\frac{\varepsilon}{2}\left(\sum_{i=1}^{3}\sigma_{i}\right)^{2}\label{eq:SO3STCH_LV1_ito}
\end{align}
Define $\bar{\sigma}=\sum_{i=1}^{3}\sigma_{i}$. Substitute $\dot{\hat{b}}$,
$\dot{\hat{\sigma}}$, and $W$ from \eqref{eq:SO3STCH_best_ito},
\eqref{eq:SO3STCH_sest_ito}, and \eqref{eq:SO3STCH_Wcorr_ito}, respectively,
in \eqref{eq:SO3STCH_LV1_ito}. Also, $||\tilde{R}||_{I}=\left\Vert \tilde{\rho}\right\Vert ^{2}/\left(1+\left\Vert \tilde{\rho}\right\Vert ^{2}\right)$
and $\boldsymbol{\Upsilon}_{a}(\tilde{R})=2\tilde{\rho}/\left(1+\left\Vert \tilde{\rho}\right\Vert ^{2}\right)$
as defined in \eqref{eq:SO3STCH_TR2} and \eqref{eq:SO3STCH_VEX_Pa},
respectively. Hence, \eqref{eq:SO3STCH_LV1_ito} yields

\begin{align}
\mathcal{L}V\leq & -4\left(\frac{8k_{2}-1}{8}\bar{\sigma}+\frac{32k_{1}-9}{32\varepsilon}\right)\frac{\left\Vert \tilde{\rho}\right\Vert ^{4}}{\left(1+\left\Vert \tilde{\rho}\right\Vert ^{2}\right)^{3}}\nonumber \\
& -\frac{4k_{1}}{\varepsilon}\frac{\left\Vert \tilde{\rho}\right\Vert ^{4}}{\left(1+\left\Vert \tilde{\rho}\right\Vert ^{2}\right)^{2}}-\frac{\left(1+3\left\Vert \tilde{\rho}\right\Vert ^{2}\right)\left\Vert \tilde{\rho}\right\Vert ^{2}\tilde{\rho}^{\top}\mathcal{Q}^{2}\tilde{\rho}}{2\left(1+\left\Vert \tilde{\rho}\right\Vert ^{2}\right)^{3}}\nonumber \\
& +k_{b}\tilde{b}^{\top}\hat{b}+k_{\sigma}\tilde{\sigma}^{\top}\hat{\sigma}+\frac{\varepsilon}{2}\bar{\sigma}^{2}\label{eq:SO3STCH_LV2_ito}
\end{align}
from \eqref{eq:SO3STCH_LV2_ito} $k_{b}\tilde{b}^{\top}\hat{b}=-k_{b}||\tilde{b}||^{2}+k_{b}\tilde{b}^{\top}b$
and $k_{\sigma}\tilde{\sigma}^{\top}\hat{\sigma}=-k_{\sigma}\left\Vert \tilde{\sigma}\right\Vert ^{2}+k_{\sigma}\tilde{\sigma}^{\top}\sigma$.
Combining this result with Young's inequality yields 
\begin{align*}
k_{b}\tilde{b}^{\top}b & \leq\frac{k_{b}}{2}||\tilde{b}||^{2}+\frac{k_{b}}{2}\left\Vert b\right\Vert ^{2}\\
k_{\sigma}\tilde{\sigma}^{\top}\sigma & \leq\frac{k_{\sigma}}{2}\left\Vert \tilde{\sigma}\right\Vert ^{2}+\frac{k_{\sigma}}{2}\left\Vert \sigma\right\Vert ^{2}
\end{align*}
thereby, the differential operator in \eqref{eq:SO3STCH_LV2_ito}
results in 
\begin{align}
\mathcal{L}V\leq & -4\left(\frac{8k_{2}-1}{8}\bar{\sigma}+\frac{32k_{1}-9}{32\varepsilon}\right)\frac{\left\Vert \tilde{\rho}\right\Vert ^{4}}{\left(1+\left\Vert \tilde{\rho}\right\Vert ^{2}\right)^{3}}\nonumber \\
& -\frac{\left(1+3\left\Vert \tilde{\rho}\right\Vert ^{2}\right)\left\Vert \tilde{\rho}\right\Vert ^{2}\tilde{\rho}^{\top}\mathcal{Q}^{2}\tilde{\rho}}{2\left(1+\left\Vert \tilde{\rho}\right\Vert ^{2}\right)^{3}}\nonumber \\
& -\frac{4k_{1}}{\varepsilon}\frac{\left\Vert \tilde{\rho}\right\Vert ^{4}}{\left(1+\left\Vert \tilde{\rho}\right\Vert ^{2}\right)^{2}}-\frac{k_{b}}{2}||\tilde{b}||^{2}-\frac{k_{\sigma}}{2}\left\Vert \tilde{\sigma}\right\Vert ^{2}\nonumber \\
& +\frac{k_{b}}{2}\left\Vert b\right\Vert ^{2}+\frac{k_{\sigma}}{2}\left\Vert \sigma\right\Vert ^{2}+\frac{\varepsilon}{2}\bar{\sigma}^{2}\label{eq:SO3STCH_LV3_ito}
\end{align}
such that \eqref{eq:SO3STCH_LV3_ito} in $\mathbb{SO}\left(3\right)$
form is equivalent to 
\begin{align}
\mathcal{L}V\leq & -\left(\frac{8k_{2}-1}{8}\bar{\sigma}+\frac{32k_{1}-9}{32\varepsilon}\right)||\tilde{R}||_{I}\left\Vert \boldsymbol{\Upsilon}_{a}(\tilde{R})\right\Vert ^{2}\nonumber \\
& -\left(\frac{1}{8}+\frac{3}{8}\frac{||\tilde{R}||_{I}}{1-||\tilde{R}||_{I}}\right)||\tilde{R}||_{I}\boldsymbol{\Upsilon}_{a}(\tilde{R})^{\top}\mathcal{Q}^{2}\boldsymbol{\Upsilon}_{a}(\tilde{R})\nonumber \\
& -\frac{4k_{1}}{\varepsilon}||\tilde{R}||_{I}^{2}-\frac{k_{b}}{2}||\tilde{b}||^{2}-\frac{k_{\sigma}}{2}\left\Vert \tilde{\sigma}\right\Vert ^{2}\nonumber \\
& +\frac{k_{b}}{2}\left\Vert b\right\Vert ^{2}+\frac{k_{\sigma}}{2}\left\Vert \sigma\right\Vert ^{2}+\frac{\varepsilon}{2}\bar{\sigma}^{2}\label{eq:SO3STCH_LV5_ito}
\end{align}
Setting $\gamma_{1}\geq1$, $\gamma_{2}\geq1$, $k_{1}\geq\frac{9}{32}$,
$k_{2}\geq\frac{1}{8}$, $k_{b}>0$, $k_{\sigma}>0$, and the positive
constant $\varepsilon$ sufficiently small with $\mathcal{Q}^{2}:\mathbb{R}_{+}\rightarrow\mathbb{R}^{3\times3}$
being bounded, the operator $\mathcal{L}V$ in \eqref{eq:SO3STCH_LV3_ito}
becomes similar to \eqref{eq:SO3STCH_dVfunction_Lyap} in Lemma \ref{lem:SO3STCH_1}.
Define $c_{2}=\frac{k_{b}}{2}\left\Vert b\right\Vert ^{2}+\frac{1}{2}\left(k_{\sigma}+\varepsilon\right)\bar{\sigma}^{2}$
which is governed by the unknown constant bias $b$ and the the upper
bound of covariance $\sigma$. Let $\tilde{X}=\left[\frac{\left\Vert \tilde{\rho}\right\Vert ^{2}}{1+\left\Vert \tilde{\rho}\right\Vert ^{2}},\frac{1}{\sqrt{2\gamma_{1}}}\tilde{b}^{\top},\frac{1}{\sqrt{2\gamma_{2}}}\tilde{\sigma}^{\top}\right]^{\top}\in\mathbb{R}^{7}$
and 
\[
\mathcal{H}=\left[\begin{array}{ccc}
4k_{1}/\varepsilon & \underline{\mathbf{0}}_{3}^{\top} & \underline{\mathbf{0}}_{3}^{\top}\\
\underline{\mathbf{0}}_{3} & \gamma_{1}k_{b}\mathbf{I}_{3} & \mathbf{0}_{3\times3}\\
\underline{\mathbf{0}}_{3} & \mathbf{0}_{3\times3} & \gamma_{2}k_{\sigma}\mathbf{I}_{3}
\end{array}\right]\in\mathbb{R}^{7\times7}
\]
Hence, the differential operator in \eqref{eq:SO3STCH_LV3_ito} can
be expressed as 
\begin{align}
\mathcal{L}V\leq & -4\left(\frac{8k_{2}-1}{8}\bar{\sigma}+\frac{32k_{1}-9}{32\varepsilon}\right)\frac{\left\Vert \tilde{\rho}\right\Vert ^{4}}{\left(1+\left\Vert \tilde{\rho}\right\Vert ^{2}\right)^{3}}\nonumber \\
& -\frac{\left(1+3\left\Vert \tilde{\rho}\right\Vert ^{2}\right)\left\Vert \tilde{\rho}\right\Vert ^{2}\tilde{\rho}^{\top}\mathcal{Q}^{2}\tilde{\rho}}{2\left(1+\left\Vert \tilde{\rho}\right\Vert ^{2}\right)^{3}}-\tilde{X}^{\top}\mathcal{H}\tilde{X}+c_{2}\label{eq:SO3STCH_LV4_ito}
\end{align}
or more simply 
\begin{equation}
\mathcal{L}V\leq-h\left(\left\Vert \tilde{\rho}\right\Vert \right)-\underline{\lambda}\left(\mathcal{H}\right)V+c_{2}\label{eq:SO3STCH_LV_Final_ito}
\end{equation}
such that $h\left(\cdot\right)$ is a class $\mathcal{K}$ function
which includes the first two components in \eqref{eq:SO3STCH_LV4_ito},
and $\underline{\lambda}\left(\cdot\right)$ denotes the minimum eigenvalue
of a matrix. Based on \eqref{eq:SO3STCH_LV_Final_ito}, one easily
obtains 
\begin{equation}
\frac{d\left(\mathbb{E}\left[V\right]\right)}{dt}=\mathbb{E}\left[\mathcal{L}V\right]\leq-\underline{\lambda}\left(\mathcal{H}\right)\mathbb{E}\left[V\right]+c_{2}\label{eq:SO3STCH_dV_Exp_ito}
\end{equation}
Consider ${\rm K}=\mathbb{E}\left[V\left(t\right)\right]$; thus $\frac{d\left(\mathbb{E}\left[V\right]\right)}{dt}\leq0$
for $\underline{\lambda}\left(\mathcal{H}\right)>\frac{c_{2}}{{\rm K}}$.
Hence, $V\leq{\rm K}$ is an invariant set and for $\mathbb{E}\left[V\left(0\right)\right]\leq{\rm K}$
there is $\mathbb{E}\left[V\left(t\right)\right]\leq{\rm K}\forall t>0$.
Based on Lemma \ref{lem:SO3STCH_1}, the inequality in \eqref{eq:SO3STCH_dV_Exp_ito}
holds for $V\left(0\right)\leq{\rm K}$ and for all $t>0$ such that
\begin{equation}
0\leq\mathbb{E}\left[V\left(t\right)\right]\leq V\left(0\right){\rm exp}\left(-\underline{\lambda}\left(\mathcal{H}\right)t\right)+\frac{c_{2}}{\underline{\lambda}\left(\mathcal{H}\right)},\,\forall t\geq0\label{eq:SO3STCH_V_Exp_ito}
\end{equation}
The above-mentioned inequality implies that $\mathbb{E}\left[V\left(t\right)\right]$
is eventually bounded by $c_{2}/\underline{\lambda}\left(\mathcal{H}\right)$
indicating that $\tilde{X}$ is SGUUB in the mean square. Let us define
$\tilde{Y}=\left[\tilde{\rho}^{\top},\tilde{b}^{\top},\tilde{\sigma}^{\top}\right]^{\top}\in\mathbb{R}^{9}$.
Since $\tilde{X}$ is SGUUB, $\tilde{Y}$ is SGUUB in the mean square.
For a rotation matrix $R\in\mathbb{SO}\left(3\right)$, let us define
$\mathcal{U}_{0}\subseteq\mathbb{SO}\left(3\right)\times\mathbb{R}^{3}\times\mathbb{R}^{3}$
as $\mathcal{U}_{0}=\left\{ \left.\left(\tilde{R}_{0},\tilde{b}_{0},\tilde{\sigma}_{0}\right)\right|{\rm Tr}\left\{ \tilde{R}_{0}\right\} =-1,\tilde{b}_{0}=\underline{\mathbf{0}}_{3},\tilde{\sigma}_{0}=\underline{\mathbf{0}}_{3}\right\} $.
The set $\mathcal{U}_{0}$ is forward invariant and unstable for the
dynamic system in \eqref{eq:SO3STCH_R_dynam}. From almost any initial
condition such that $\tilde{R}_{0}\notin\mathcal{U}_{0}$ or, equivalently,
$\tilde{\rho}_{0}\in\mathbb{R}^{3}$, the trajectory of $\tilde{X}$
is SGUUB in the mean square. 

\subsection{Nonlinear Stochastic Attitude Filter in Stratonovich Sense}

Stochastic differential equations can be defined and solved in the
sense of Ito integral \cite{ito1984lectures}. Alternatively, Stratonovich
integral \cite{stratonovich1967topics} can be employed for solving
stochastic differential equations. The common feature between Stratonovich
and Ito integral is that if the associated function multiplied by
$d\beta$ is continuous and Lipschitz, the mean square limit exists.
The Ito integral is defined for functional on $\left\{ \beta\left(\tau\right),\tau\leq t\right\} $
which is more natural but does not obey the chain rule. Conversely,
Stratonovich is a well-defined Riemann integral for the sampled function,
it has a continuous partial derivative with respect to $\beta$, it
obeys the chain rule and it is more convenient for colored noise \cite{stratonovich1967topics,jazwinski2007stochastic}.
Hence, the Stratonovich integral is defined for explicit functions
of $\beta$. In case when angular velocity measurements are contaminated
with a wide-band of random colored noise process, the solution of
\eqref{eq:SO3STCH_Stoch_non1} for $\rho\left(t_{0}\right)=0$ is
defined by 
\begin{equation}
\rho\left(t\right)=\int_{t_{0}}^{t}f\left(\rho\left(\tau\right),b\left(\tau\right)\right)d\tau+\int_{t_{0}}^{t}g\left(\rho\left(\tau\right)\right)\mathcal{Q}d\beta\label{eq:SO3STCH_Strat_solution}
\end{equation}
according to subsection \ref{subsec:SO3STCH_Nonlinear-Ito}, the expected
value of \eqref{eq:SO3STCH_Strat_solution} is 
\[
\mathbb{E}\left[\rho\right]\neq\int_{t_{0}}^{t}\mathbb{E}\left[f\left(\rho\left(\tau\right),b\left(\tau\right)\right)\right]d\tau
\]
Thus, Stratonovich introduced the Wong-Zakai correction factor which
can help in designing an adaptive estimate for the covariance component.
Let us assume that the attitude dynamic in \eqref{eq:SO3STCH_Stoch_Ito}
was defined in the sense of Stratonovich \cite{stratonovich1967topics},
hence, its equivalent Ito \cite{ito1984lectures,khasminskii1980stochastic,jazwinski2007stochastic}
can be defined by 
\begin{align}
\left[d\rho\right]_{i}= & \left[f\left(\rho,b\right)\right]_{i}dt+\sum_{k=1}^{3}\sum_{j=1}^{3}\frac{\mathcal{Q}_{j,j}^{2}}{2}g_{kj}\left(\rho\right)\frac{\partial g_{ij}\left(\rho\right)}{\partial\rho_{k}}dt\nonumber \\
& +\left[g\left(\rho\right)\mathcal{Q}d\beta\right]_{i}\label{eq:SO3STCH_Stoch_strat_i}
\end{align}
where both $f\left(\rho,b\right)$ and $g\left(\rho\right)$ are defined
in \eqref{eq:SO3STCH_Stoch_Ito}, $i,j,k=1,2,3$ denote $i$th, $j$th
and $k$th element components of the associate vector or matrix. The
term $\sum_{k=1}^{3}\sum_{j=1}^{3}\frac{\mathcal{Q}_{j,j}^{2}}{2}g_{kj}\left(\rho\right)\frac{\partial g_{ij}\left(\rho\right)}{\partial\rho_{k}}$
denotes the Wong-Zakai correction factor of stochastic differential
equations (SDEs) in the sense of Ito's representations \cite{wong1965convergence}.
Let $\boldsymbol{\mathcal{W}}_{i}\left(\rho\right)=\sum_{k=1}^{3}\sum_{j=1}^{3}\frac{\mathcal{Q}_{j,j}^{2}}{2}g_{kj}\left(\rho\right)\frac{\partial g_{ij}\left(\rho\right)}{\partial\rho_{k}}$,
accordingly, one can find that for $i=1$ 
\begin{align*}
& \sum_{k=1}^{3}\sum_{j=1}^{3}\frac{\mathcal{Q}_{j,j}^{2}}{2}g_{kj}\left(\rho\right)\frac{\partial g_{ij}\left(\rho\right)}{\partial\rho_{k}}=\frac{1}{4}\left(\left(1+\rho_{1}^{2}\right)\rho_{1}\mathcal{Q}_{1,1}^{2}+\right.\\
& \hspace{7em}\left.\left(\rho_{1}\rho_{2}-\rho_{3}\right)\rho_{2}\mathcal{Q}_{2,2}^{2}+\left(\rho_{2}+\rho_{1}\rho_{3}\right)\rho_{3}\mathcal{Q}_{3,3}^{2}\right)
\end{align*}
Hence, $\boldsymbol{\mathcal{W}}\left(\rho\right)$ for $i=1,2,3$
is 
\begin{align}
\boldsymbol{\mathcal{W}}\left(\rho\right) & =\frac{1}{4}\left(\mathbf{I}_{3}+\left[\rho\right]_{\times}+\rho\rho^{\top}\right)\mathcal{Q}^{2}\rho\label{eq:SO3STCH_WONG-Zaki}
\end{align}
Manipulating equations \eqref{eq:SO3STCH_Stoch_strat_i} and \eqref{eq:SO3STCH_WONG-Zaki},
the stochastic dynamics of the Rodriguez vector can be expressed as
\begin{align}
d\rho= & \mathcal{F}\left(\rho,b\right)dt+g\left(\rho\right)\mathcal{Q}d\beta\label{eq:SO3STCH_Stoch_strat}
\end{align}
where $g\left(\rho\right):=-\frac{1}{2}\left(\mathbf{I}_{3}+\left[\rho\right]_{\times}+\rho\rho^{\top}\right)$
and $\mathcal{F}\left(\rho,b\right):=-g\left(\rho\right)\left(\Omega_{m}-b\right)+\boldsymbol{\mathcal{W}}\left(\rho\right)$.
Define the error in attitude estimate similar to \eqref{eq:SO3STCH_R_error}.
Also, assume that the elements of covariance matrix $\mathcal{Q}^{2}$
are upper bounded by $\sigma$ as given in \eqref{eq:SO3STCH_g_factor}
such that the bound of $\sigma$ is unknown with nonnegative elements.

Consider the following nonlinear stochastic filter design 
\begin{align}
\dot{\hat{R}}= & \hat{R}\left[\Omega_{m}-\hat{b}-\frac{1}{2}\frac{{\rm diag}\left(\boldsymbol{\Upsilon}_{a}(\tilde{R})\right)}{1-||\tilde{R}||_{I}}\hat{\sigma}-W\right]_{\times}\label{eq:SO3STCH_dRest_Strat}\\
\dot{\hat{b}}= & \gamma_{1}||\tilde{R}||_{I}\boldsymbol{\Upsilon}_{a}(\tilde{R})-\gamma_{1}k_{b}\hat{b},\hspace{1em}\hat{b}\left(0\right)=\hat{b}_{0}\label{eq:SO3STCH_best_Strat}\\
\mathcal{\dot{\hat{\sigma}}}= & \gamma_{2}||\tilde{R}||_{I}\left(k_{1}\mathcal{D}_{\Upsilon}^{\top}+\frac{1}{2}\frac{{\rm diag}\left(\boldsymbol{\Upsilon}_{a}(\tilde{R})\right)}{1-||\tilde{R}||_{I}}\right)\boldsymbol{\Upsilon}_{a}(\tilde{R})\nonumber \\
& -\gamma_{2}k_{\sigma}\hat{\sigma},\hspace{1em}\hat{\sigma}\left(0\right)=\hat{\sigma}_{0}\label{eq:SO3STCH_sest_strat}\\
W= & \frac{k_{1}}{\varepsilon}\frac{2-||\tilde{R}||_{I}}{1-||\tilde{R}||_{I}}\boldsymbol{\Upsilon}_{a}(\tilde{R})+k_{2}\mathcal{D}_{\Upsilon}\hat{\sigma}\label{eq:SO3STCH_Wcorr_Strat}
\end{align}
where $\hat{R}\left(0\right)=\hat{R}_{0}$, $\Omega_{m}$ is the angular
velocity measurement as defined in \eqref{eq:SO3STCH_Angular}, $\hat{b}$
and $\hat{\sigma}$ are estimates of the unknown parameters $b$ and
$\sigma$, respectively, $\tilde{R}=R_{y}^{\top}\hat{R}$ with $R_{y}$
being the reconstructed attitude, $\boldsymbol{\Upsilon}_{a}(\tilde{R})=\mathbf{vex}\left(\boldsymbol{\mathcal{P}}_{a}(\tilde{R})\right)$
was given in \eqref{eq:SO3STCH_VEX_Pa}, $||\tilde{R}||_{I}$ is the
Euclidean distance of $\tilde{R}$, and $\mathcal{D}_{\Upsilon}=\left[\boldsymbol{\Upsilon}_{a}(\tilde{R}),\boldsymbol{\Upsilon}_{a}(\tilde{R}),\boldsymbol{\Upsilon}_{a}(\tilde{R})\right]$.
$\gamma_{1}$ and $\gamma_{2}$ are positive adaptation gains, $\varepsilon>0$
is a small constant, while $k_{b}$, $k_{\sigma}$, $k_{1}$ and $k_{2}$
are positive constants. Quaternion representation of Stratonovich's filter is presented in \nameref{sec:SO3STCH_AppendixB}. 
\begin{thm}
	\label{thm:SO3STCH_2} Consider the rotation kinematics in \eqref{eq:SO3STCH_R_Dyn1}
	with angular velocity measurements and given vector measurements in
	\eqref{eq:SO3STCH_Angular} and \eqref{eq:SO3STCH_Vector_norm}, respectively,
	being coupled with the observer in \eqref{eq:SO3STCH_dRest_Strat},
	\eqref{eq:SO3STCH_best_Strat}, \eqref{eq:SO3STCH_sest_strat} and
	\eqref{eq:SO3STCH_Wcorr_Strat}. Assume that two or more body-frame
	non-collinear vectors are available for measurements. Then, for angular
	velocity measurements contaminated with noise $\left(\omega\neq\underline{\mathbf{0}}_{3}\right)$, $\left[\tilde{\rho}^{\top},\tilde{b}^{\top},\tilde{\sigma}^{\top}\right]^{\top}$is
	semi-globally uniformly ultimately bounded in mean square. Moreover,
	the observer errors can be made sufficiently small by choosing the
	appropriate design parameters.
\end{thm}
\textbf{Proof: }Let the error in vector $b$ and $\sigma$ be defined
as in \eqref{eq:SO3STCH_b_tilde} and \eqref{eq:SO3STCH_sigma_tilde_strat},
respectively. Hence, the derivative of \eqref{eq:SO3STCH_R_error}
in incremental form can be obtained from \eqref{eq:SO3STCH_Stoch_Ito}
and \eqref{eq:SO3STCH_dRest_Strat} by 
\begin{align}
d\tilde{R}= & \tilde{R}\left[\Omega-\tilde{R}^{\top}\Omega+\tilde{b}-\frac{1}{2}{\rm diag}\left(\tilde{\rho}\right)\hat{\sigma}-W\right]_{\times}dt\nonumber \\
& +\tilde{R}\left[\mathcal{Q}d\beta\right]_{\times}\label{eq:SO3STCH_dR_tilde_strat}
\end{align}
Assume that the Rodriguez error vector dynamic of \eqref{eq:SO3STCH_dR_tilde_strat}
is defined in the sense of Stratonovich. The extraction of Rodriguez
error vector dynamics in view of the transformation of \eqref{eq:SO3STCH_R_dynam}
into \eqref{eq:SO3STCH_Stoch_strat} can be expressed from \eqref{eq:SO3STCH_dR_tilde_strat}
to \eqref{eq:SO3STCH_StochErr_Strat} in Ito's representation \cite{stratonovich1967topics}
as 
\begin{equation}
d\tilde{\rho}=\tilde{\mathcal{F}}dt+\tilde{g}\mathcal{Q}d\beta\label{eq:SO3STCH_StochErr_Strat}
\end{equation}
where $\tilde{\rho}$ is Rodriguez error vector associated with $\tilde{R}$
with $\tilde{g}=\frac{1}{2}\left(\mathbf{I}_{3}+\left[\tilde{\rho}\right]_{\times}+\tilde{\rho}\tilde{\rho}^{\top}\right)$,
$\tilde{\mathcal{F}}=\tilde{g}\left(\Omega-\tilde{R}^{\top}\Omega+\tilde{b}-\frac{1}{2}{\rm diag}\left(\tilde{\rho}\right)\hat{\sigma}-W\right)+\boldsymbol{\mathcal{W}}\left(\tilde{\rho}\right)$
and $\boldsymbol{\mathcal{W}}\left(\tilde{\rho}\right)=\frac{1}{4}\left(\mathbf{I}_{3}+\left[\tilde{\rho}\right]_{\times}+\tilde{\rho}\tilde{\rho}^{\top}\right)\mathcal{Q}^{2}\tilde{\rho}$.
Consider the following potential function 
\begin{align}
V\left(\tilde{\rho},\tilde{b},\tilde{\sigma}\right)= & \left(\frac{\left\Vert \tilde{\rho}\right\Vert ^{2}}{1+\left\Vert \tilde{\rho}\right\Vert ^{2}}\right)^{2}+\frac{1}{2\gamma_{1}}\tilde{b}^{\top}\tilde{b}+\frac{1}{2\gamma_{2}}\tilde{\sigma}^{\top}\tilde{\sigma}\label{eq:SO3STCH_LyapV_strat}
\end{align}
For $V:=V\left(\tilde{\rho},\tilde{b},\tilde{\sigma}\right)$, the
differential operator $\mathcal{L}V$ in Definition \ref{def:SO3STCH_2}
for the dynamic system in \eqref{eq:SO3STCH_StochErr_Strat} can be
written as 
\begin{equation}
\mathcal{L}V=V_{\tilde{\rho}}^{\top}\tilde{\mathcal{F}}+\frac{1}{2}{\rm Tr}\left\{ \tilde{g}^{\top}V_{\tilde{\rho}\tilde{\rho}}\tilde{g}\mathcal{Q}^{2}\right\} -\frac{1}{\gamma_{1}}\tilde{b}^{\top}\dot{\hat{b}}-\frac{1}{\gamma_{2}}\tilde{\sigma}^{\top}\dot{\hat{\sigma}}\label{eq:SO3STCH_LV_strat}
\end{equation}
The first and the second partial derivatives of \eqref{eq:SO3STCH_LyapV_strat}
with respect to $\tilde{\rho}$ are similar to \eqref{eq:SO3STCH_LyapVv_ito}
and \eqref{eq:SO3STCH_LyapVvv_ito}, respectively. The first part
of differential operator $\mathcal{L}V$ in \eqref{eq:SO3STCH_LV_strat}
can be evaluated by 
\begin{align}
V_{\mathcal{\tilde{\rho}}}^{\top}\tilde{\mathcal{F}}= & 2\frac{\left\Vert \tilde{\rho}\right\Vert ^{2}}{\left(1+\left\Vert \tilde{\rho}\right\Vert ^{2}\right)^{2}}\tilde{\rho}^{\top}\left(\tilde{b}-\frac{1}{2}{\rm diag}\left(\tilde{\rho}\right)\hat{\sigma}+\frac{1}{2}\mathcal{Q}^{2}\tilde{\rho}-W\right)\nonumber \\
\leq & 2\frac{\left\Vert \tilde{\rho}\right\Vert ^{2}}{\left(1+\left\Vert \tilde{\rho}\right\Vert ^{2}\right)^{2}}\tilde{\rho}^{\top}\left(\tilde{b}+\frac{1}{2}{\rm diag}\left(\tilde{\rho}\right)\tilde{\sigma}-W\right)\label{eq:SO3STCH_Vdot_strat}
\end{align}
where $\frac{1}{2}V_{\mathcal{\tilde{\rho}}}^{\top}\left(\mathbf{I}_{3}+\left[\tilde{\rho}\right]_{\times}+\tilde{\rho}\tilde{\rho}^{\top}\right)\left(\Omega-\tilde{R}^{\top}\Omega\right)=0$.
The component ${\rm Tr}\left\{ \tilde{g}^{\top}V_{\tilde{\rho}\tilde{\rho}}\tilde{g}\mathcal{Q}^{2}\right\} $
is similar to \eqref{eq:SO3STCH_LyapTrace_ito}. Let $\bar{q}=\left[\mathcal{Q}_{1,1},\mathcal{Q}_{2,2},\mathcal{Q}_{3,3}\right]^{\top}$
and $\sigma$ be similar to \eqref{eq:SO3STCH_g_factor}. The operator
$\mathcal{L}V$ in \eqref{eq:SO3STCH_LyapV_strat} becomes 
\begin{align*}
& \mathcal{L}V\leq2\frac{\left\Vert \tilde{\rho}\right\Vert ^{2}\tilde{\rho}^{\top}\left(\tilde{b}-W+\frac{1}{2}{\rm diag}\left(\tilde{\rho}\right)\tilde{\sigma}\right)}{\left(1+\left\Vert \tilde{\rho}\right\Vert ^{2}\right)^{2}}+\frac{{\rm Tr}\left\{ \left\Vert \tilde{\rho}\right\Vert ^{2}\mathcal{Q}^{2}\right\} }{2\left(1+\left\Vert \tilde{\rho}\right\Vert ^{2}\right)^{2}}\\
& +\frac{{\rm Tr}\left\{ \left(2-\left\Vert \tilde{\rho}\right\Vert ^{2}-3\left\Vert \tilde{\rho}\right\Vert ^{4}\right)\tilde{\rho}\tilde{\rho}^{\top}\mathcal{Q}^{2}\right\} }{2\left(1+\left\Vert \tilde{\rho}\right\Vert ^{2}\right)^{3}}-\frac{1}{\gamma_{1}}\tilde{b}^{\top}\dot{\hat{b}}-\frac{1}{\gamma_{2}}\tilde{\sigma}^{\top}\dot{\hat{\sigma}}
\end{align*}
Since $\left\Vert \bar{q}\right\Vert ^{2}={\rm Tr}\left\{ \mathcal{Q}^{2}\right\} $
and ${\rm Tr}\left\{ \tilde{\rho}\tilde{\rho}^{\top}\mathcal{Q}^{2}\right\} \leq\left\Vert \tilde{\rho}\right\Vert ^{2}\left\Vert \bar{q}\right\Vert ^{2}$,
we obtain 
\begin{align}
\mathcal{L}V\leq & 2\frac{\left\Vert \tilde{\rho}\right\Vert ^{2}\tilde{\rho}^{\top}\left(\tilde{b}-W+\frac{1}{2}{\rm diag}\left(\tilde{\rho}\right)\tilde{\sigma}\right)}{\left(1+\left\Vert \tilde{\rho}\right\Vert ^{2}\right)^{2}}\nonumber \\
& +\frac{\left\Vert \tilde{\rho}\right\Vert ^{4}{\rm Tr}\left\{ \mathcal{Q}^{2}\right\} +3\left\Vert \tilde{\rho}\right\Vert ^{2}\left\Vert \bar{q}\right\Vert ^{2}}{2\left(1+\left\Vert \tilde{\rho}\right\Vert ^{2}\right)^{3}}-\frac{1}{\gamma_{1}}\tilde{b}^{\top}\dot{\hat{b}}\nonumber \\
& -\frac{1}{\gamma_{2}}\tilde{\sigma}^{\top}\dot{\hat{\sigma}}-\frac{\left\Vert \tilde{\rho}\right\Vert ^{2}\left(1+3\left\Vert \tilde{\rho}\right\Vert ^{2}\right)\tilde{\rho}^{\top}\mathcal{Q}^{2}\tilde{\rho}}{2\left(1+\left\Vert \tilde{\rho}\right\Vert ^{2}\right)^{3}}\label{eq:SO3STCH_LV1_strat}
\end{align}
From the last result and taking into consideration the inequality
in \eqref{eq:SO3STCH_Inq_Lyap}, according to Lemma \ref{lem:SO3STCH_2},
and \eqref{eq:SO3STCH_g_factor}, equation \eqref{eq:SO3STCH_LV1_strat}
becomes 
\begin{align}
\mathcal{L}V\leq & 2\frac{\left\Vert \tilde{\rho}\right\Vert ^{2}\tilde{\rho}^{\top}\left(\tilde{b}-W+\frac{1}{2}{\rm diag}\left(\tilde{\rho}\right)\tilde{\sigma}\right)}{\left(1+\left\Vert \tilde{\rho}\right\Vert ^{2}\right)^{2}}\nonumber \\
& +2\frac{\left\Vert \tilde{\rho}\right\Vert ^{2}\tilde{\rho}^{\top}\left(\frac{1}{4}\mathcal{D}_{\tilde{\rho}}\sigma+\frac{9}{16\varepsilon}\tilde{\rho}\right)}{\left(1+\left\Vert \tilde{\rho}\right\Vert ^{2}\right)^{3}}-\frac{1}{\gamma_{1}}\tilde{b}^{\top}\dot{\hat{b}}-\frac{1}{\gamma_{2}}\tilde{\sigma}^{\top}\dot{\hat{\sigma}}\nonumber \\
& -\frac{\left\Vert \tilde{\rho}\right\Vert ^{2}\left(1+3\left\Vert \tilde{\rho}\right\Vert ^{2}\right)\tilde{\rho}^{\top}\mathcal{Q}^{2}\tilde{\rho}}{2\left(1+\left\Vert \tilde{\rho}\right\Vert ^{2}\right)^{3}}+\frac{\varepsilon}{2}\left(\sum_{i=1}^{3}\sigma_{i}\right)^{2}\label{eq:SO3STCH_LV2_strat}
\end{align}
with $\mathcal{D}_{\tilde{\rho}}=\left[\tilde{\rho},\tilde{\rho},\tilde{\rho}\right]$.
From \eqref{eq:SO3STCH_LV2_strat}, we have $\tilde{\rho}^{\top}\mathcal{D}_{\tilde{\rho}}\sigma=\left(\sum_{i=1}^{3}\sigma_{i}\right)\left\Vert \tilde{\rho}\right\Vert ^{2}$.
Let us define $\bar{\sigma}=\sum_{i=1}^{3}\sigma_{i}$. Substitute
for the differential operators $\dot{\hat{b}}$, $\dot{\hat{\sigma}}$
and the correction factor $W$ from \eqref{eq:SO3STCH_best_Strat},
\eqref{eq:SO3STCH_sest_strat} and \eqref{eq:SO3STCH_Wcorr_Strat},
respectively, with $||\tilde{R}||_{I}=\left\Vert \tilde{\rho}\right\Vert ^{2}/\left(1+\left\Vert \tilde{\rho}\right\Vert ^{2}\right)$
and $\boldsymbol{\Upsilon}_{a}(\tilde{R})=2\tilde{\rho}/\left(1+\left\Vert \tilde{\rho}\right\Vert ^{2}\right)$.
Hence, the result in \eqref{eq:SO3STCH_LV2_strat} is equivalent to
\begin{align}
\mathcal{L}V\leq & -4\left(\frac{8k_{2}-1}{8}\bar{\sigma}+\frac{32k_{1}-9}{32\varepsilon}\right)\frac{\left\Vert \tilde{\rho}\right\Vert ^{4}}{\left(1+\left\Vert \tilde{\rho}\right\Vert ^{2}\right)^{3}}\nonumber \\
& -\frac{\left(1+3\left\Vert \tilde{\rho}\right\Vert ^{2}\right)\left\Vert \tilde{\rho}\right\Vert ^{2}\tilde{\rho}^{\top}\mathcal{Q}^{2}\tilde{\rho}}{2\left(1+\left\Vert \tilde{\rho}\right\Vert ^{2}\right)^{3}}\nonumber \\
& -\frac{4k_{1}}{\varepsilon}\frac{\left\Vert \tilde{\rho}\right\Vert ^{4}}{\left(1+\left\Vert \tilde{\rho}\right\Vert ^{2}\right)^{2}}-k_{b}||\tilde{b}||^{2}-k_{\sigma}\left\Vert \tilde{\sigma}\right\Vert ^{2}\nonumber \\
& +k_{b}\tilde{b}^{\top}b+k_{\sigma}\tilde{\sigma}^{\top}\sigma+\frac{\varepsilon}{2}\bar{\sigma}^{2}\label{eq:SO3STCH_LV3_strat}
\end{align}
applying Young’s inequality, one has 
\begin{align*}
k_{b}\tilde{b}^{\top}b & \leq\frac{k_{b}}{2}||\tilde{b}||^{2}+\frac{k_{b}}{2}\left\Vert b\right\Vert ^{2}\\
k_{\sigma}\tilde{\sigma}^{\top}\sigma & \leq\frac{k_{\sigma}}{2}\left\Vert \tilde{\sigma}\right\Vert ^{2}+\frac{k_{\sigma}}{2}\bar{\sigma}^{2}
\end{align*}
Consequently, \eqref{eq:SO3STCH_LV3_strat} becomes{\small{} 
	\begin{align}
	\mathcal{L}V\leq & -4\left(\frac{8k_{2}-1}{8}\bar{\sigma}+\frac{32k_{1}-9}{32\varepsilon}\right)\frac{\left\Vert \tilde{\rho}\right\Vert ^{4}}{\left(1+\left\Vert \tilde{\rho}\right\Vert ^{2}\right)^{3}}\nonumber \\
	& -\frac{\left(1+3\left\Vert \tilde{\rho}\right\Vert ^{2}\right)\left\Vert \tilde{\rho}\right\Vert ^{2}\tilde{\rho}^{\top}\mathcal{Q}^{2}\tilde{\rho}}{2\left(1+\left\Vert \tilde{\rho}\right\Vert ^{2}\right)^{3}}-\frac{4k_{1}}{\varepsilon}\frac{\left\Vert \tilde{\rho}\right\Vert ^{4}}{\left(1+\left\Vert \tilde{\rho}\right\Vert ^{2}\right)^{2}}-\frac{k_{b}}{2}||\tilde{b}||^{2}\nonumber \\
	& -\frac{k_{\sigma}}{2}\left\Vert \tilde{\sigma}\right\Vert ^{2}+\frac{k_{b}}{2}\left\Vert b\right\Vert ^{2}+\frac{1}{2}\left(k_{\sigma}+\varepsilon\right)\bar{\sigma}^{2}\label{eq:SO3STCH_LV_Final_strat}
	\end{align}
}In other words, \eqref{eq:SO3STCH_LV_Final_strat} in $\mathbb{SO}\left(3\right)$
form is equivalent to{\small{} 
	\begin{align}
	\mathcal{L}V\leq & -\left(\frac{1}{8}+\frac{3}{8}\frac{||\tilde{R}||_{I}}{1-||\tilde{R}||_{I}}\right)||\tilde{R}||_{I}\boldsymbol{\Upsilon}_{a}(\tilde{R})^{\top}\mathcal{Q}^{2}\boldsymbol{\Upsilon}_{a}(\tilde{R})\nonumber \\
	& -\left(\frac{8k_{2}-1}{8}\bar{\sigma}+\frac{32k_{1}-9}{32\varepsilon}\right)||\tilde{R}||_{I}\left\Vert \boldsymbol{\Upsilon}_{a}(\tilde{R})\right\Vert ^{2}-\frac{4k_{1}}{\varepsilon}||\tilde{R}||_{I}^{2}\nonumber \\
	& -\frac{k_{b}}{2}||\tilde{b}||^{2}-\frac{k_{\sigma}}{2}\left\Vert \tilde{\sigma}\right\Vert ^{2}+\frac{k_{b}}{2}\left\Vert b\right\Vert ^{2}+\frac{1}{2}\left(k_{\sigma}+\varepsilon\right)\bar{\sigma}^{2}\label{eq:SO3STCH_LV4_strat}
	\end{align}
}Setting $\gamma_{1}\geq1$, $\gamma_{2}\geq1$, $k_{1}\geq\frac{9}{32}$,
$k_{2}\geq\frac{1}{8}$, $k_{b}>0$, $k_{\sigma}>0$, and the positive
constant $\varepsilon$ being sufficiently small, and defining $c_{2}=\frac{k_{b}}{2}\left\Vert b\right\Vert ^{2}+\frac{1}{2}\left(k_{\sigma}+\varepsilon\right)\bar{\sigma}^{2}$,
the operator $\mathcal{L}V$ in \eqref{eq:SO3STCH_LV_Final_strat}
becomes similar to (4.16) in \cite{deng2001stabilization} which is
in turn similar to \eqref{eq:SO3STCH_dVfunction_Lyap} in Lemma \ref{lem:SO3STCH_1}.
Define
\begin{align*}
\tilde{X} & =\left[\frac{\left\Vert \tilde{\rho}\right\Vert ^{2}}{1+\left\Vert \tilde{\rho}\right\Vert ^{2}},\frac{1}{\sqrt{2\gamma_{1}}}\tilde{b}^{\top},\frac{1}{\sqrt{2\gamma_{2}}}\tilde{\sigma}^{\top}\right]^{\top}\in\mathbb{R}^{7},\\
\mathcal{H} & =\left[\begin{array}{ccc}
4k_{1}/\varepsilon & \underline{\mathbf{0}}_{3}^{\top} & \underline{\mathbf{0}}_{3}^{\top}\\
\underline{\mathbf{0}}_{3} & \gamma_{1}k_{b}\mathbf{I}_{3} & \mathbf{0}_{3\times3}\\
\underline{\mathbf{0}}_{3} & \mathbf{0}_{3\times3} & \gamma_{2}k_{\sigma}\mathbf{I}_{3}
\end{array}\right]\in\mathbb{R}^{7\times7}
\end{align*}
Thereby, the differential operator in \eqref{eq:SO3STCH_LV_Final_strat}
is{\small{} 
	\begin{align}
	\mathcal{L}V\leq & -\left(\frac{8k_{2}-1}{2}\bar{\sigma}+\frac{32k_{1}-9}{8\varepsilon}\right)\frac{\left\Vert \tilde{\rho}\right\Vert ^{4}}{\left(1+\left\Vert \tilde{\rho}\right\Vert ^{2}\right)^{3}}\nonumber \\
	& -\frac{\left(1+3\left\Vert \tilde{\rho}\right\Vert ^{2}\right)\left\Vert \tilde{\rho}\right\Vert ^{2}\tilde{\rho}^{\top}\mathcal{Q}^{2}\tilde{\rho}}{2\left(1+\left\Vert \tilde{\rho}\right\Vert ^{2}\right)^{3}}-\tilde{X}^{\top}\mathcal{H}\tilde{X}+c_{2}\nonumber \\
	\leq & -h\left(\left\Vert \tilde{\rho}\right\Vert \right)-\underline{\lambda}\left(\mathcal{H}\right)V+c_{2}\label{eq:SO3STCH_LV5_strat}
	\end{align}
}such that $h\left(\cdot\right)$ is a class $\mathcal{K}$ function
which includes the first two components in \eqref{eq:SO3STCH_LV5_strat}.
Based on \eqref{eq:SO3STCH_LV5_strat}, one easily obtains 
\begin{equation}
\frac{d\left(\mathbb{E}\left[V\right]\right)}{dt}=\mathbb{E}\left[\mathcal{L}V\right]\leq-\underline{\lambda}\left(\mathcal{H}\right)\mathbb{E}\left[V\right]+c_{2}\label{eq:SO3STCH_dV_Exp_strat}
\end{equation}
Let ${\rm K}=\mathbb{E}\left[V\left(t\right)\right]$; then $\frac{d\left(\mathbb{E}\left[V\right]\right)}{dt}\leq0$
for $\underline{\lambda}\left(\mathcal{H}\right)>\frac{c_{2}}{{\rm K}}$.
Thereby, $V\left(t\right)\leq{\rm K}$ is an invariant set and for
$\mathbb{E}\left[V\left(0\right)\right]\leq{\rm K}$ it follows that
$\mathbb{E}\left[V\left(t\right)\right]\leq{\rm K}\forall t>0$. Accordingly,
the inequality in \eqref{eq:SO3STCH_dV_Exp_strat} holds for $V\left(0\right)\leq{\rm K}$
and for all $t>0$ which means that 
\begin{equation}
0\leq\mathbb{E}\left[V\left(t\right)\right]\leq V\left(0\right){\rm exp}\left(-\underline{\lambda}\left(\mathcal{H}\right)t\right)+\frac{c_{2}}{\underline{\lambda}\left(\mathcal{H}\right)},\,\forall t\geq0\label{eq:SO3STCH_V_Exp_strat}
\end{equation}
The above inequality entails that $\mathbb{E}\left[V\left(t\right)\right]$
is eventually bounded by $c_{2}/\underline{\lambda}\left(\mathcal{H}\right)$
which implies that $\tilde{X}$ is SGUUB in the mean square. For a
rotation matrix $R\in\mathbb{SO}\left(3\right)$, define $\mathcal{U}_{0}\subseteq\mathbb{SO}\left(3\right)\times\mathbb{R}^{3}\times\mathbb{R}^{3}$
as $\mathcal{U}_{0}=\left\{ \left.\left(\tilde{R}_{0},\tilde{b}_{0},\tilde{\sigma}_{0}\right)\right|{\rm Tr}\left\{ \tilde{R}_{0}\right\} =-1,\tilde{b}_{0}=\underline{\mathbf{0}}_{3},\tilde{\sigma}_{0}=\underline{\mathbf{0}}_{3}\right\} $.
The set $\mathcal{U}_{0}$ is forward invariant and unstable for the
dynamic system in \eqref{eq:SO3STCH_R_dynam}. Therefore, for almost
any initial condition such that $\tilde{R}_{0}\notin\mathcal{U}_{0}$
or, equivalently, for any $\tilde{\rho}_{0}\in\mathbb{R}^{3}$, $\tilde{X}$
is SGUUB in the mean square as in Definition \ref{def:SO3STCH_1}.

Since, $\mathcal{Q}^{2}:\mathbb{R}_{+}\rightarrow\mathbb{R}^{3\times3}$
is bounded, we have $d\left(\mathbb{E}\left[V\right]\right)/dt<0$
for $V>c_{2}/\underline{\lambda}\left(\mathcal{H}\right)$. Considering
Lemma \ref{lem:SO3STCH_1} and the design parameters of the stochastic
observer in Theorem \ref{thm:SO3STCH_1} or \ref{thm:SO3STCH_2} and
combining them with prior knowledge about the covariance upper bound,
allows to make the error signal smaller if the design parameters are
chosen appropriately. 

\subsection{Stochastic Attitude Filters: Ito vs Stratonovich \label{sec:SO3STCH_subsec_comp}}

In this work, the selection of potential functions in \eqref{eq:SO3STCH_LyapV_ito}
and \eqref{eq:SO3STCH_LyapV_strat} contributes to attenuating and
controlling the noise level associated with angular velocity measurements.
Also, the selection of potential functions in \eqref{eq:SO3STCH_LyapV_ito}
and \eqref{eq:SO3STCH_LyapV_strat} produced results analogous to
those \eqref{eq:SO3STCH_LV_Final_ito} and \eqref{eq:SO3STCH_LV5_strat},
respectively. This similarity in potential function selection and
final results is critical as it guarantees fair comparison between
the two proposed stochastic filters. The proposed stochastic filters
are able to correct the attitude allowing the user to reduce the noise
level associated with angular velocity measurements through $\underline{\lambda}\left(\mathcal{H}\right)$
by setting the values of $\varepsilon$, $k_{1}$, $k_{b}$, $k_{\sigma}$,
$\gamma_{1}$ and $\gamma_{2}$ appropriately. Nonlinear deterministic
attitude filters lack this advantage.\\
The main features of the nonlinear stochastic attitude filter in
the sense of Ito can be listed as 
\begin{enumerate}
	\item[1)] The filter requires less computational power in comparison with the
	Stratonovich's filter.
	\item[2)] No prior information about the covariance matrix $\mathcal{Q}^{2}$
	is required. 
	\item[3)] This filter is applicable to white noise. 
\end{enumerate}
Whereas, the main characteristics of the nonlinear stochastic attitude
filter in the sense of Stratonovich are 
\begin{enumerate}
	\item[1)] The filter demands more computational power in comparison with the
	Ito's filter. 
	\item[2)] No prior information about the covariance matrix $\mathcal{Q}^{2}$
	is required. 
	\item[3)] The filter is applicable for white as well as colored noise. 
\end{enumerate}

\section{Simulations \label{sec:SO3STCH_Simulations}}

This section presents the performance and comparison among the two
proposed nonlinear stochastic filters on $\mathbb{SO}\left(3\right)$.
The first nonlinear stochastic filter is driven in the sense of Ito
and the second one considers Stratonovich. Consider the orientation
matrix $R$ obtained from attitude dynamics in equation \eqref{eq:SO3STCH_R_dynam}
with the following angular velocity input signal 
\[
\Omega=\left[\begin{array}{c}
{\rm sin}\left(0.7t\right)\\
0.7{\rm sin}\left(0.5t+\pi\right)\\
0.5{\rm sin}\left(0.3t+\frac{\pi}{3}\right)
\end{array}\right]\left({\rm rad/sec}\right)
\]
while the initial attitude is $R\left(0\right)=\mathbf{I}_{3}$. Let
the true angular velocity $\left(\Omega\right)$ be contaminated with
a wide-band of random noise process with zero mean and standard deviation
(STD) be equal to $0.5\left({\rm rad/sec}\right)$ such that $\Omega_{m}=\Omega+b+\omega$
with $b=0.1\left[1,-1,1\right]^{\top}$, $\omega=0.5n\left(t\right)$,
where $t$ denotes real time, $n\left(t\right)={\rm randn}\left(3,1\right)$
where ${\rm randn}\left(3,1\right)$ is a ${\rm MATLAB}{}^{\circledR}$
command, which refers to a normally distributed random vector at each
time instant. Let non-collinear inertial-frame vectors be given as
${\rm v}_{1}^{\mathcal{I}}=\frac{1}{\sqrt{3}}\left[1,-1,1\right]^{\top}$
and ${\rm v}_{2}^{\mathcal{I}}=\left[0,0,1\right]^{\top}$, while
body-frame vectors ${\rm v}_{1}^{\mathcal{B}}$ and ${\rm v}_{2}^{\mathcal{B}}$
are obtained by ${\rm v}_{i}^{\mathcal{B}}=R^{\top}{\rm v}_{i}^{\mathcal{I}}+{\rm b}_{i}^{\mathcal{B}}+\omega_{i}^{\mathcal{B}}$
for $i=1,2$. Also, suppose that an additional noise vector $\omega_{i}^{\mathcal{B}}$
with zero mean and STD of $0.15$ corrupted the body-frame vector
measurements and bias components ${\rm b}_{1}^{\mathcal{B}}=0.1\left[-1,1,0.5\right]^{\top}$
and ${\rm b}_{2}^{\mathcal{B}}=0.1\left[0,0,1\right]^{\top}$. The
third vector of inertial-frame and body-frame is extracted by ${\rm v}_{3}^{\mathcal{I}}={\rm v}_{1}^{\mathcal{I}}\times{\rm v}_{2}^{\mathcal{I}}$
and ${\rm v}_{3}^{\mathcal{B}}={\rm v}_{1}^{\mathcal{B}}\times{\rm v}_{2}^{\mathcal{B}}$
and followed by normalization of the three vectors at each time instant
according to \eqref{eq:SO3STCH_Vector_norm}. From vectorial measurements,
the corrupted reconstructed attitude $R_{y}$ is obtained by SVD \cite{markley1988attitude}
with $\tilde{R}=R_{y}^{\top}\hat{R}$, see \nameref{sec:SO3STCH_AppendixA}.
The total simulation time is 15 seconds.

For a very large initial attitude error, the initial rotation of attitude
estimate is given according to angle-axis parameterization in \eqref{eq:SO3STCH_att_ang}
by $\hat{R}\left(0\right)=\mathcal{R}_{\alpha}\left(\alpha,u/\left\Vert u\right\Vert \right)$
with $\alpha=179.9\left({\rm deg}\right)$ and $u=\left[1,5,3\right]^{\top}$
being very close to the unstable equilibria such that $||\tilde{R}\left(0\right)||_{I}\approx0.99999$.
The initial conditions are{\footnotesize{}
	\[
	R\left(0\right)=\left[\begin{array}{ccc}
	1 & 0 & 0\\
	0 & 1 & 0\\
	0 & 0 & 1
	\end{array}\right],\hspace{1em}\hat{R}\left(0\right)=\left[\begin{array}{ccc}
	-0.9429 & 0.2848 & 0.1729\\
	0.2866 & 0.4286 & 0.8568\\
	0.1700 & 0.8574 & -0.4857
	\end{array}\right]
	\]
}Initial estimates for both filters are $\hat{b}\left(0\right)=\left[0,0,0\right]^{\top}$
and $\hat{\sigma}\left(0\right)=\left[0,0,0\right]^{\top}$. The same
notation is used in derivations of both nonlinear stochastic filters.
The design parameters were chosen as $\gamma_{1}=1$, $\gamma_{2}=1$,
$k_{b}=0.5$, $k_{\sigma}=0.5$, $k_{1}=0.5$, $k_{2}=0.5$ and $\varepsilon=0.5$.
Additionally, the following color notation is used: green color demonstrates
the true value, red illustrates the performance of Ito's filter and
blue represents the performance of Stratonovich stochastic filter.
Also, magenta refers to a measured value.

\begin{figure}[h!]
	\centering{}\includegraphics[scale=0.26]{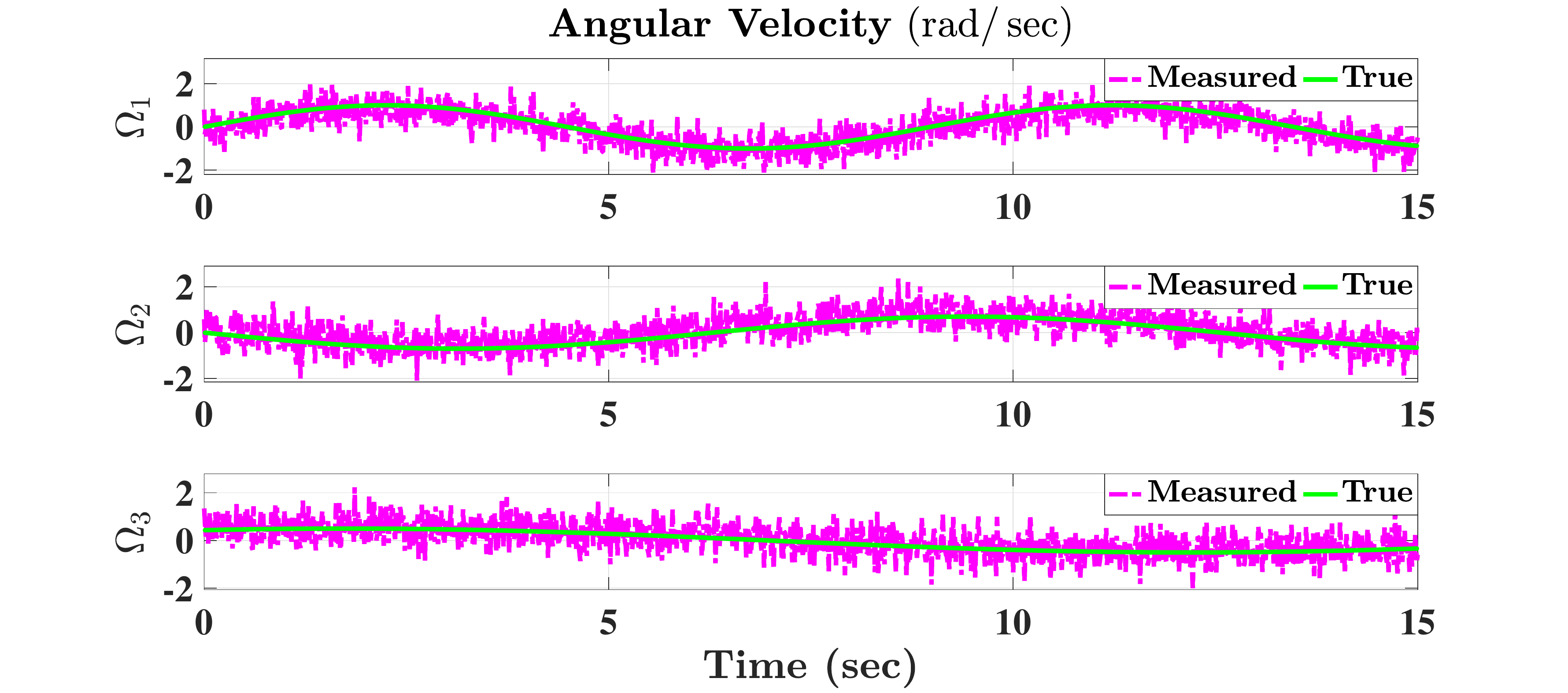}\caption{True and measured angular velocities.}
	\label{fig:SO3STCH_1} 
\end{figure}

\begin{figure}[h!]
	\centering{}\includegraphics[scale=0.26]{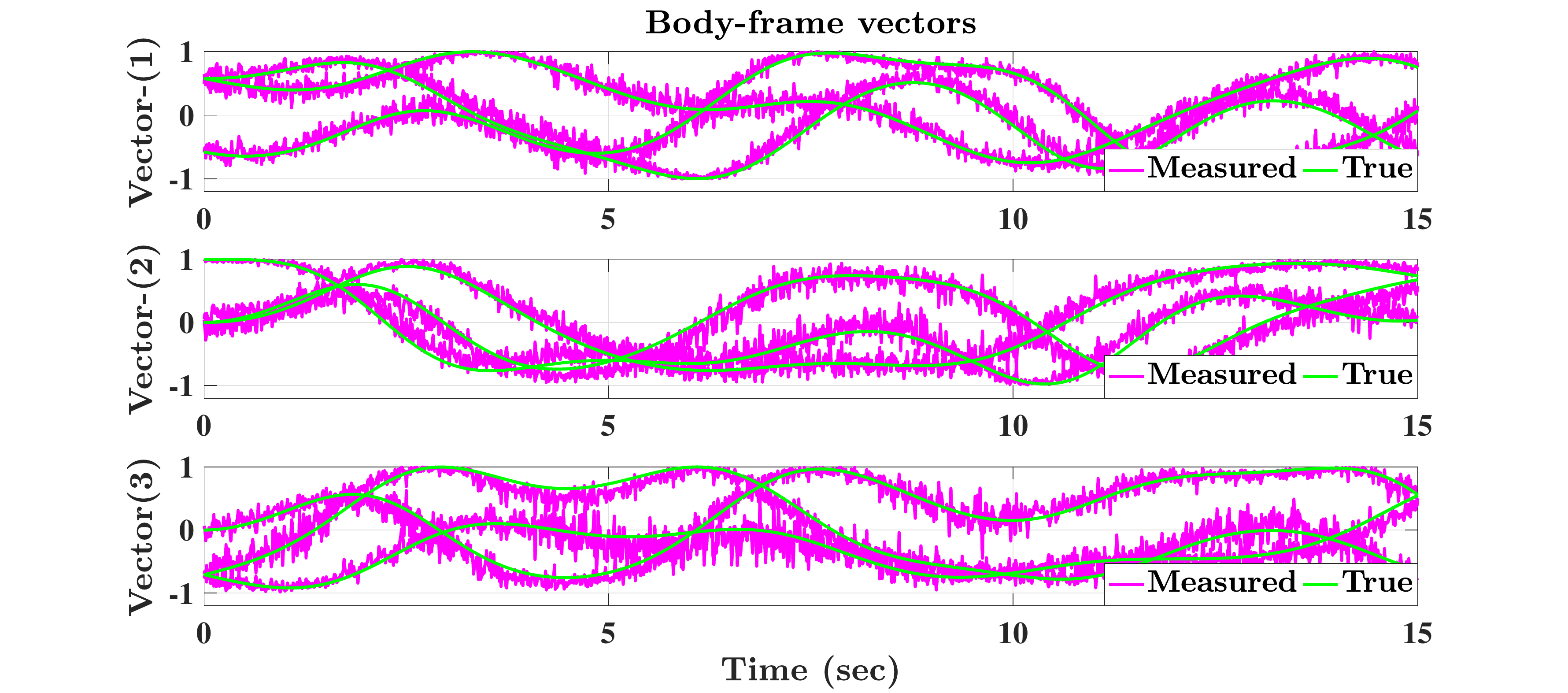}\caption{True and measured body-frame vectorial measurements.}
	\label{fig:SO3STCH_2} 
\end{figure}

\begin{figure}[h!]
	\centering{}\includegraphics[scale=0.26]{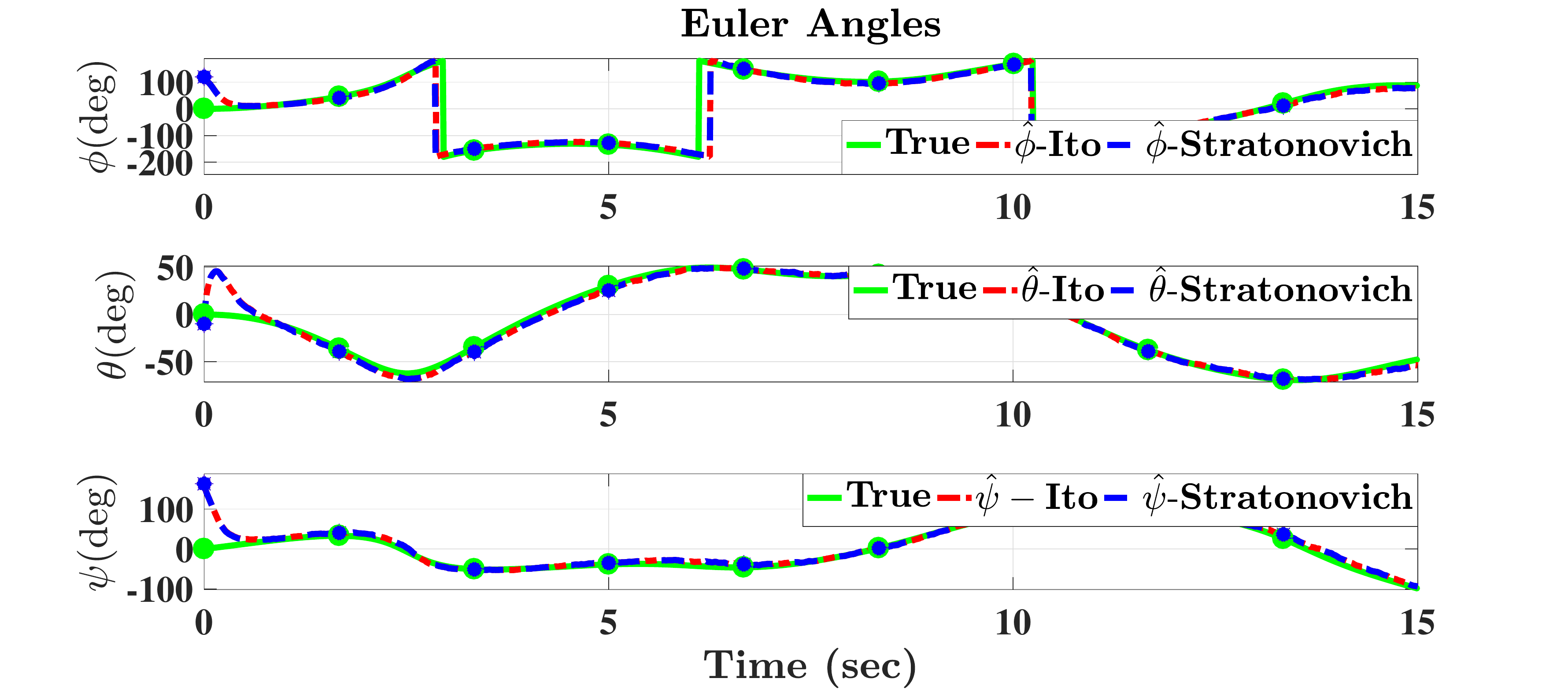}\caption{Tracking performance of Euler angles.}
	\label{fig:SO3STCH_3} 
\end{figure}

\begin{figure}[h!]
	\centering{}\includegraphics[scale=0.26]{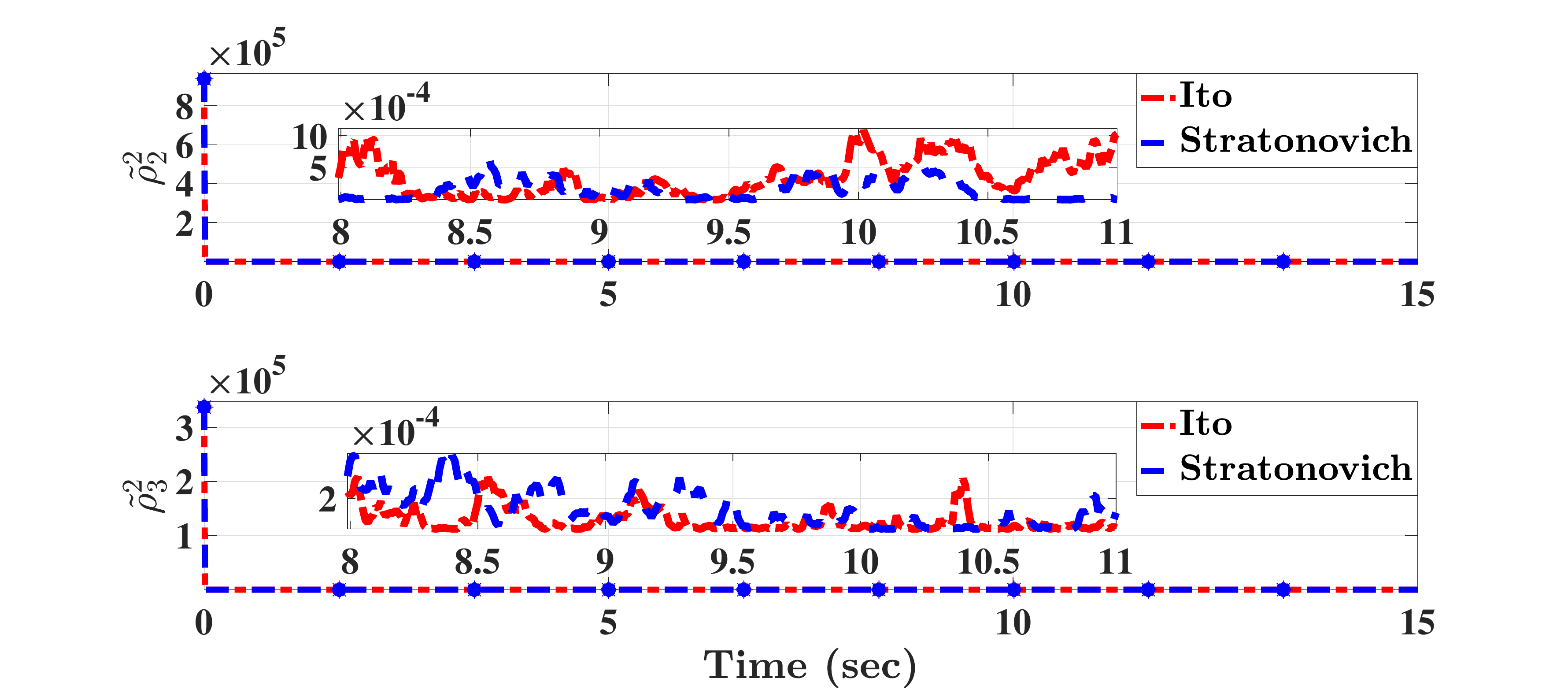}\caption{Rodriguez vector square error $\tilde{\rho}^{2}$.}
	\label{fig:SO3STCH_4} 
\end{figure}

\begin{figure}[h]
	\centering{}\includegraphics[scale=0.26]{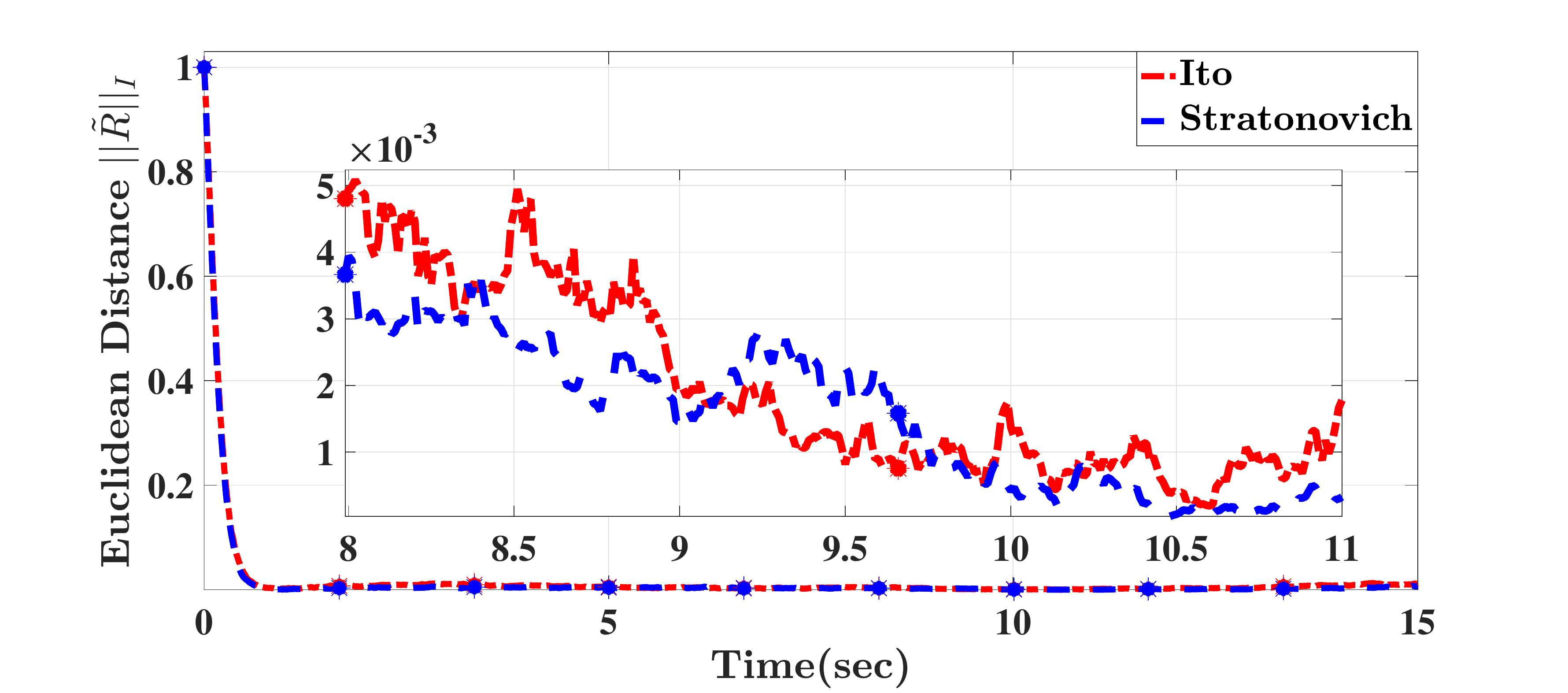}\caption{Tracking performance of normalized Euclidean distance error $||\tilde{R}||_{I}$.}
	\label{fig:SO3STCH_5} 
\end{figure}

The true angular velocity $\left(\Omega\right)$ and the high values
of noise and bias components introduced through the measurement process
of $\Omega_{m}$ plotted against time are depicted in Fig. \ref{fig:SO3STCH_1}.
Also, Fig. \ref{fig:SO3STCH_2} presents the true body-frame vectors
and their uncertain measurements. Fig. \ref{fig:SO3STCH_3} shows
the tracked Euler angles $\left(\phi,\theta,\psi\right)$ of Ito and Stratonovich
stochastic attitude filters relative to true angles plotted against
time. Fig. \ref{fig:SO3STCH_3} presents impressive tracking performance
of the proposed stochastic filters. The mapping from $\mathbb{SO}\left(3\right)$
implies that $\tilde{\rho}\rightarrow\infty$ as $||\tilde{R}||_{I}\rightarrow1$.
Accordingly, Fig. \ref{fig:SO3STCH_4} demonstrates the convergence
of the square error of Rodriguez vector $\tilde{\rho}^{2}$ from large
error initialization to a very small value close to zero. Fig. \ref{fig:SO3STCH_5}
confirms all the previous discussion using normalized Euclidean distance
$||\tilde{R}||_{I}=\frac{1}{4}{\rm Tr}\left\{ \mathbf{I}_{3}-R^{\top}\hat{R}\right\} $
which shows remarkable stable and fast convergence to very small neighborhood
of the origin. However, Ito stochastic filter is characterized by
higher oscillatory performance compared to Stratonovich stochastic
filter.

\vspace{8pt}
\vspace{8pt}
\vspace{8pt}

\begin{figure}[h]
	\centering{}\includegraphics[scale=0.26]{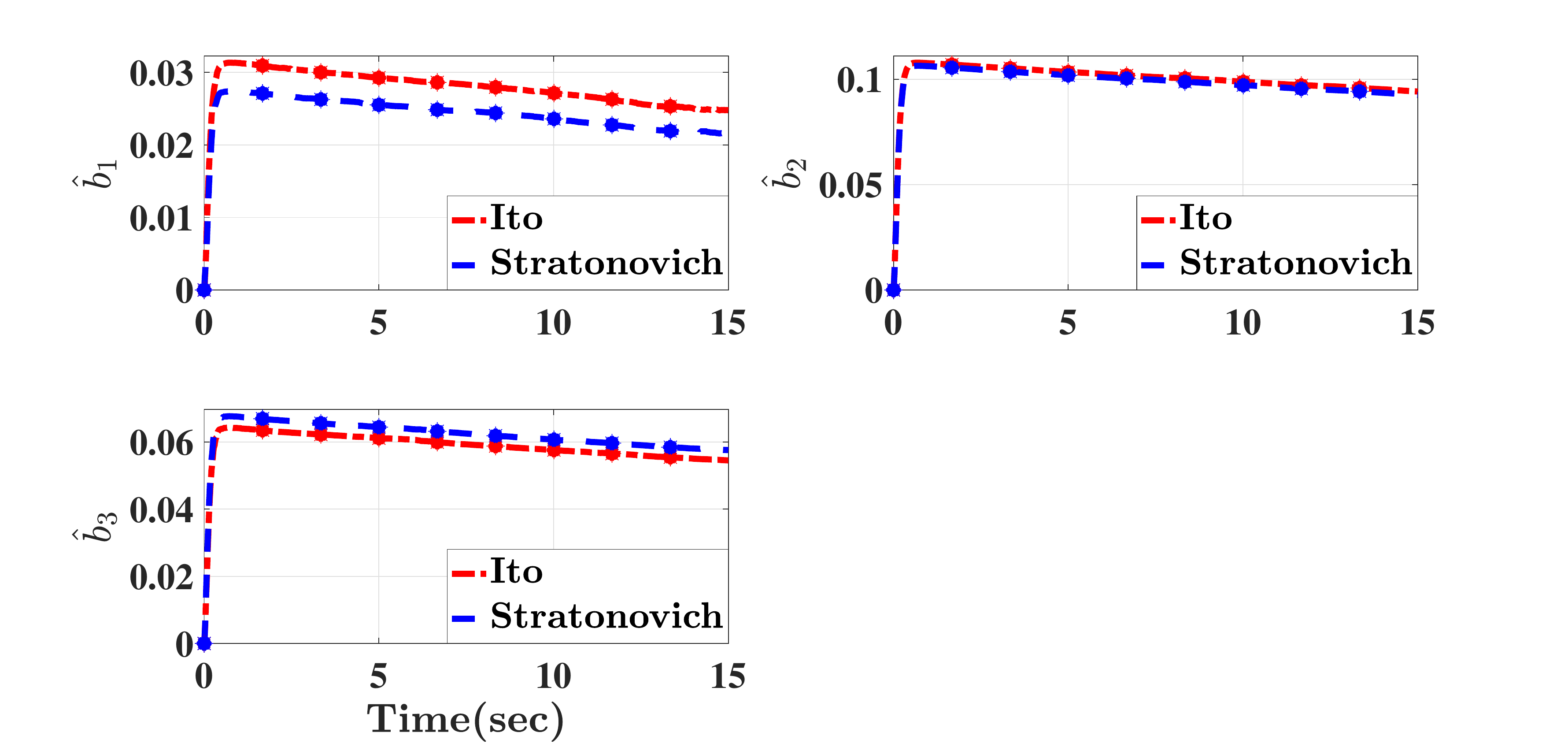}\caption{Estimates of stochastic attitude filters ($\hat{b}$).}
	\label{fig:SO3STCH_6} 
\end{figure}

\begin{figure}[h]
	\centering{}\includegraphics[scale=0.26]{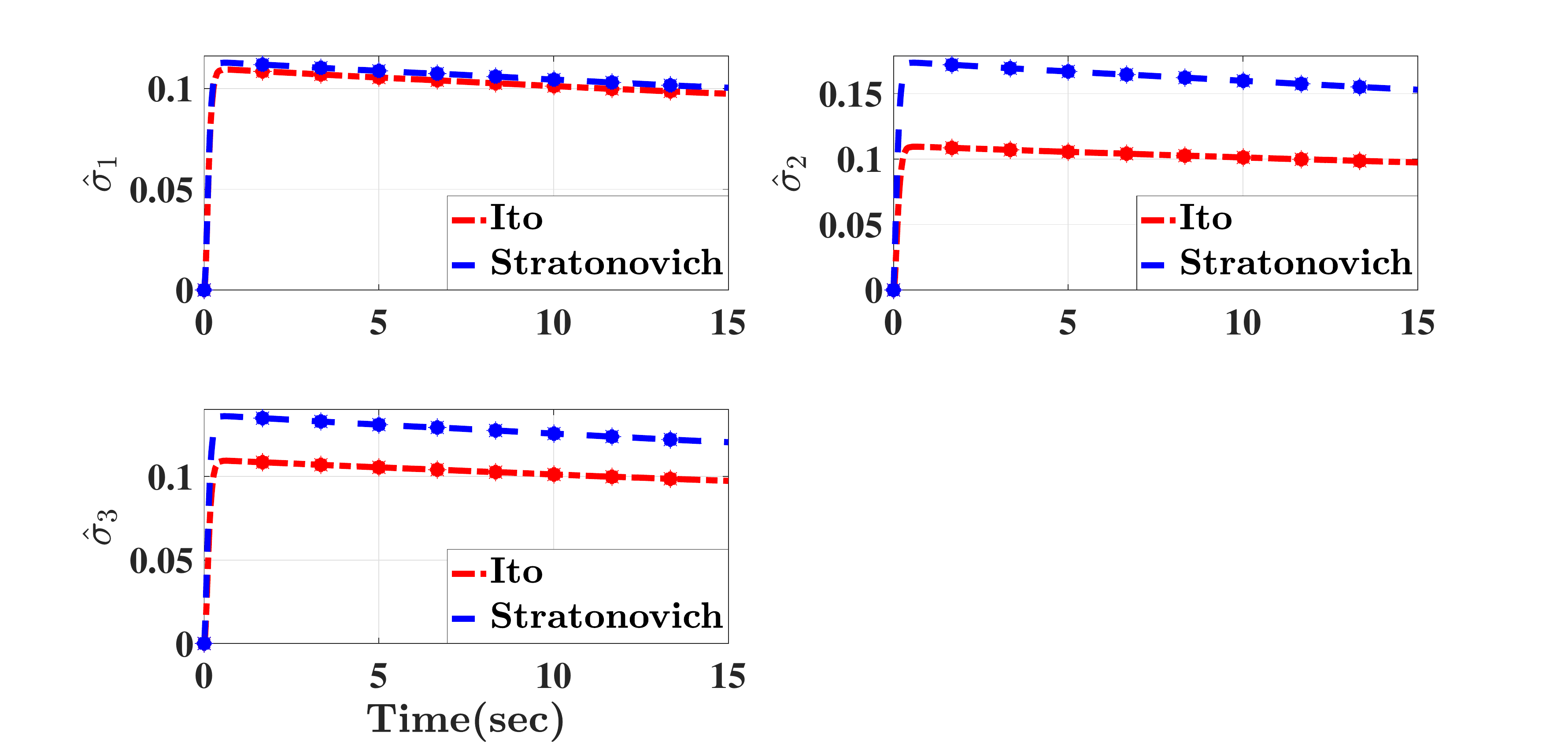}\caption{Estimates of stochastic attitude filters ($\hat{\sigma}$).}
	\label{fig:SO3STCH_7} 
\end{figure}

To further compare the steady-state performance of the proposed filters
in terms of normalized Euclidean distance of the error ($||\tilde{R}||_{I}$),
Table \ref{tab:SO3STCH_1} summarizes statistical details of the mean
and the STD of $||\tilde{R}||_{I}$. Both filters showed very small
mean error of $||\tilde{R}||_{I}$ with $||\tilde{R}||_{I}$ being
regulated to close neighborhood of the origin however, Stratonovich's
filter showed a remarkable less mean errors and STD in comparison
with Ito's filter. Numerical results included in Table \ref{tab:SO3STCH_1}
proves that the proposed nonlinear stochastic filters are robust as
illustrated in Fig. \ref{fig:SO3STCH_3}, \ref{fig:SO3STCH_4}, and
\ref{fig:SO3STCH_5}.
\vspace{10pt}
\begin{table}[H]
	\caption{Statistical analysis of $||\tilde{R}||_{I}$ of the two proposed filter.}
	
	\begin{centering}
		\begin{tabular}{c|>{\centering}p{2.5cm}|>{\centering}p{2.5cm}}
			\hline 
			\multicolumn{3}{c}{Output data of $||\tilde{R}||_{I}$ over the period (1-15 sec)}\tabularnewline
			\hline 
			\hline 
			Filter  & Ito  & Stratonovich\tabularnewline
			\hline 
			Mean  & $4.1\times10^{-3}$  & $2.8\times10^{-3}$\tabularnewline
			\hline 
			STD  & $3\times10^{-3}$  & $1.6\times10^{-3}$\tabularnewline
			\hline 
		\end{tabular}
		\par\end{centering}
	\label{tab:SO3STCH_1} 
\end{table}
\vspace{10pt}
Finally, Fig. \ref{fig:SO3STCH_6} and \ref{fig:SO3STCH_7} illustrate
the estimates of the stochastic filters plotted against time. It can
be concluded from Fig. \ref{fig:SO3STCH_6} and \ref{fig:SO3STCH_7}
that the estimates of the proposed filter are stable and smooth.

\newpage
Results show effectiveness and robustness of the two stochastic filters
against bias and noise components contaminating angular velocity measurements,
as well as uncertainty in vectorial measurements and large initial
error. Stochastic filters have proven to be able to correct their
attitude in a small amount of time requiring no prior information
about the covariance matrix $\mathcal{Q}^{2}$ in order to obtain
impressive estimation performance. The main advantage of Stratonovich
stochastic filter, as mentioned in Subsection \ref{sec:SO3STCH_subsec_comp},
is that the filter is applicable to white as well as colored noise.
In addition, it had smaller mean square error and STD to Ito's filter
as given in Table \ref{tab:SO3STCH_1}. Nonetheless, Ito stochastic
filter requires less computational power.

\section{Conclusion \label{sec:SO3STCH_Conclusion}}

Deterministic filters neglect the noise associated with the angular
velocity measurements in filter derivation. This can be clearly noticed
in the selection of the potential function. However, an alternate
potential function which has not been considered in the literature
is able to significantly attenuate the effects of noise in angular
velocities to lower levels. As such, this paper reformulated the attitude
problem to stochastic sense through Rodriguez vector parameterization.
Two different nonlinear stochastic attitude filters on the Special
Orthogonal Group 3 ($\mathbb{SO}\left(3\right)$) have been proposed. The first filter
is developed in the sense of Ito and the second filter is driven in
the sense of Stratonovich. The resulting estimators have proven to
have fast convergence properties in the presence of high levels of
noise in angular velocity and vectorial measurements.



\section*{Acknowledgment}

The authors would like to thank University of Western Ontario for
the funding that made this research possible. Also, the authors would
like to thank \textbf{Maria Shaposhnikova} for proofreading the article.

\section*{Appendix A \label{sec:SO3STCH_AppendixA}}
\begin{center}
	\textbf{\large{}{}{}{}{}{}{}{}{}{}{}{}An Overview on SVD
		in }{\large{}{}{}{}{}{}{}{}{}\cite{markley1988attitude} {}{}
	} 
	\par\end{center}

Let $R\in\mathbb{SO}\left(3\right)$ be the true attitude. The attitude
can be reconstructed through a set of vectors given in \eqref{eq:SO3STCH_Vector_norm}.
Let $s_{i}$ be the confidence level of measurement $i$ such that
for $n$ measurements we have $\sum_{i=1}^{n}s_{i}=1$. In that case,
the corrupted reconstructed attitude $R_{y}$ can be obtained by 
\[
\begin{cases}
\mathcal{J}\left(R\right) & =1-\sum_{i=1}^{n}s_{i}\left(\upsilon_{i}^{\mathcal{B}}\right)^{\top}R^{\top}\upsilon_{i}^{\mathcal{I}}\\
& =1-{\rm Tr}\left\{ R^{\top}B^{\top}\right\} \\
B & =\sum_{i=1}^{n}s_{i}\upsilon_{i}^{\mathcal{B}}\left(\upsilon_{i}^{\mathcal{I}}\right)^{\top}=USV^{\top}\\
U_{+} & =U\left[\begin{array}{ccc}
1 & 0 & 0\\
0 & 1 & 0\\
0 & 0 & {\rm det}\left(U\right)
\end{array}\right]\\
V_{+} & =V\left[\begin{array}{ccc}
1 & 0 & 0\\
0 & 1 & 0\\
0 & 0 & {\rm det}\left(V\right)
\end{array}\right]\\
R_{y} & =V_{+}U_{+}^{\top}
\end{cases}
\]
For more details visit \cite{markley1988attitude}.

\section*{Appendix B \label{sec:SO3STCH_AppendixB}}
\begin{center}
	\textbf{\large{}{}{}{}{}{}{}{}{}{}{}{}Quaternion Representation}{\large{}{}{}
	} 
	\par\end{center}

\noindent Define $Q=[q_{0},q^{\top}]^{\top}\in\mathbb{S}^{3}$ as
a unit-quaternion with $q_{0}\in\mathbb{R}$ and $q\in\mathbb{R}^{3}$
such that $\mathbb{S}^{3}=\{\left.Q\in\mathbb{R}^{4}\right|||Q||=\sqrt{q_{0}^{2}+q^{\top}q}=1\}$.
$Q^{-1}=[\begin{array}{cc}
q_{0} & -q^{\top}\end{array}]^{\top}\in\mathbb{S}^{3}$ denotes the inverse of $Q$. Define $\odot$ as a quaternion product
where the quaternion multiplication of $Q_{1}=[\begin{array}{cc}
q_{01} & q_{1}^{\top}\end{array}]^{\top}\in\mathbb{S}^{3}$ and $Q_{2}=[\begin{array}{cc}
q_{02} & q_{2}^{\top}\end{array}]^{\top}\in\mathbb{S}^{3}$ is $Q_{1}\odot Q_{2}=[q_{01}q_{02}-q_{1}^{\top}q_{2},q_{01}q_{2}+q_{02}q_{1}+[q_{1}]_{\times}q_{2}]$.
The mapping from unit-quaternion ($\mathbb{S}^{3}$) to $\mathbb{SO}\left(3\right)$
is described by $\mathcal{R}_{Q}:\mathbb{S}^{3}\rightarrow\mathbb{SO}\left(3\right)$
\begin{align}
\mathcal{R}_{Q} & =(q_{0}^{2}-||q||^{2})\mathbf{I}_{3}+2qq^{\top}+2q_{0}\left[q\right]_{\times}\in\mathbb{SO}\left(3\right)\label{eq:NAV_Append_SO3}
\end{align}
The quaternion identity is described by $Q_{{\rm I}}=[1,0,0,0]^{\top}$
with $\mathcal{R}_{Q_{{\rm I}}}=\mathbf{I}_{3}$. Define the estimate
of $Q=[q_{0},q^{\top}]^{\top}\in\mathbb{S}^{3}$ as $\hat{Q}=[\hat{q}_{0},\hat{q}^{\top}]^{\top}\in\mathbb{S}^{3}$
with $\mathcal{R}_{\hat{Q}}=(\hat{q}_{0}^{2}-||\hat{q}||^{2})\mathbf{I}_{3}+2\hat{q}\hat{q}^{\top}+2\hat{q}_{0}\left[\hat{q}\right]_{\times}$,
see the map in \eqref{eq:NAV_Append_SO3}. The equivalent quaternion
representation of the Ito's filter in \eqref{eq:SO3STCH_dRest_ito},
\eqref{eq:SO3STCH_best_ito}, \eqref{eq:SO3STCH_sest_ito} and \eqref{eq:SO3STCH_Wcorr_ito}
is:
\[
\begin{cases}
\left[\begin{array}{c}
0\\
\upsilon_{i}^{\mathcal{B}}
\end{array}\right] & =Q^{-1}\odot\left[\begin{array}{c}
0\\
\upsilon_{i}^{\mathcal{I}}
\end{array}\right]\odot Q\\
Q_{y} & :\text{Reconstructed by QUEST algorithm}\\
\tilde{Q} & =[\tilde{q}_{0},\tilde{q}^{\top}]^{\top}=Q_{y}^{-1}\odot\hat{Q}\\
\mathcal{D}_{\Upsilon} & =2\tilde{q}_{0}\left[\tilde{q},\tilde{q},\tilde{q}\right]\\
\Gamma & =\Omega_{m}-\hat{b}-W\\
\dot{\hat{Q}} & =\frac{1}{2}\left[\begin{array}{cc}
0 & -\Gamma^{\top}\\
\Gamma & -\left[\Gamma\right]_{\times}
\end{array}\right]\hat{Q}\\
\dot{\hat{b}} & =2\gamma_{1}(1-\tilde{q}_{0}^{2})\tilde{q}_{0}\tilde{q}-\gamma_{1}k_{b}\hat{b}\\
\mathcal{\dot{\hat{\sigma}}} & =2k_{1}\gamma_{2}(1-\tilde{q}_{0}^{2})\mathcal{D}_{\Upsilon}^{\top}\tilde{q}_{0}\tilde{q}-\gamma_{2}k_{\sigma}\hat{\sigma}\\
W & =\frac{2k_{1}}{\varepsilon}\frac{1+\tilde{q}_{0}^{2}}{\tilde{q}_{0}}\tilde{q}+k_{2}\mathcal{D}_{\Upsilon}\hat{\sigma}
\end{cases}
\]
The equivalent quaternion representation of Stratonovich's filter
in \eqref{eq:SO3STCH_dRest_ito}, \eqref{eq:SO3STCH_best_ito}, \eqref{eq:SO3STCH_sest_ito}
and \eqref{eq:SO3STCH_Wcorr_ito} is:
\[
\begin{cases}
\left[\begin{array}{c}
0\\
\upsilon_{i}^{\mathcal{B}}
\end{array}\right] & =Q^{-1}\odot\left[\begin{array}{c}
0\\
\upsilon_{i}^{\mathcal{I}}
\end{array}\right]\odot Q\\
Q_{y} & :\text{Reconstructed by QUEST algorithm}\\
\tilde{Q} & =[\tilde{q}_{0},\tilde{q}^{\top}]^{\top}=Q_{y}^{-1}\odot\hat{Q}\\
\mathcal{D}_{\Upsilon} & =2\tilde{q}_{0}\left[\tilde{q},\tilde{q},\tilde{q}\right]\\
\Gamma & =\Omega_{m}-\hat{b}-\frac{{\rm diag}\left(\tilde{q}\right)}{\tilde{q}_{0}}\hat{\sigma}-W\\
\dot{\hat{Q}} & =\frac{1}{2}\left[\begin{array}{cc}
0 & -\Gamma^{\top}\\
\Gamma & -\left[\Gamma\right]_{\times}
\end{array}\right]\hat{Q}\\
\dot{\hat{b}} & =2\gamma_{1}(1-\tilde{q}_{0}^{2})\tilde{q}_{0}\tilde{q}-\gamma_{1}k_{b}\hat{b}\\
\mathcal{\dot{\hat{\sigma}}} & =2k_{1}\gamma_{2}\tilde{q}_{0}(1-\tilde{q}_{0}^{2})\mathcal{D}_{\Upsilon}^{\top}\tilde{q}\\
& \hspace{1em}+2\gamma_{2}(1-\tilde{q}_{0}^{2}){\rm diag}\left(\tilde{q}\right)\tilde{q}-\gamma_{2}k_{\sigma}\hat{\sigma}\\
W & =\frac{2k_{1}}{\varepsilon}\frac{1+\tilde{q}_{0}^{2}}{\tilde{q}_{0}}\tilde{q}+k_{2}\mathcal{D}_{\Upsilon}\hat{\sigma}
\end{cases}
\]

\bibliographystyle{IEEEtran}
\bibliography{bib_Stoch_SO3}

\vspace{310pt}

\section*{AUTHOR INFORMATION}
\vspace{10pt}
{\bf Hashim A. Hashim} is a Ph.D. candidate and a Teaching and Research Assistant in Robotics and Control, Department of Electrical and Computer Engineering at the University of Western Ontario, ON, Canada.\\
His current research interests include stochastic and deterministic filters on SO(3) and SE(3), control of multi-agent systems, control applications and optimization techniques.\\
\underline{Contact Information}: \href{mailto:hmoham33@uwo.ca}{hmoham33@uwo.ca}.
\vspace{50pt}

{\bf Lyndon J. Brown} received the B.Sc. degree from the U. of Waterloo, Canada in 1988 and the M.Sc. and PhD. degrees from the University of Illinois, Urbana-Champaign in 1991 and 1996, respectively. He is an associate professor in the department of electrical and computer engineering at Western University, Canada. He worked in industry for Honeywell Aerospace Canada and E.I. DuPont de Nemours.\\
His current research includes the identification and control of predictable signals, biological control systems, welding control systems and attitude estimation.
\vspace{50pt}

{\bf Kenneth McIsaac} received the B.Sc. degree from the University of Waterloo, Canada, in 1996, and the M.Sc. and Ph.D. degrees from the University of Pennsylvania, in 1998 and 2001, respectively. He is currently an Associate Professor and the Chair of Electrical and Computer Engineering with Western University, ON, Canada. \\
His current research interests include computer vision and signal processing, mostly in the context of using machine intelligence in robotics and assistive systems.

\end{document}